\documentclass[10pt,3p]{elsarticle}

\usepackage{standalone}
\usepackage{amssymb}
\usepackage{amsmath}
\usepackage[]{algorithm2e}

\usepackage{hyperref}
\usepackage{placeins} 	%zwecks \floatbarrier

\usepackage{todonotes}
 \usepackage[ComplainAboutUnexternalized,PatchTodo,Publish,Externaldir={./}]{externaltikz}  
\usepackage{relinput}
\journal{Journal of Computational Physics}

\newtheorem{Definition}{Definition}

\newproof{proof}{Proof}
\newcommand{\R}{\mathbb{R}}

\newcommand\MyFloatBarrier{\FloatBarrier}
\renewcommand\MyFloatBarrier{}
\begin{document}
\newcommand{\CK}{Cauchy-Kowalevsky }
\newcommand{\KC}{Kowalevsky-Cauchy }
\newcommand{\figref}[1]{Figure \ref{#1}}
\newcommand{\tabref}[1]{Table \ref{#1}}
%schreibweise prüfen
\newcommand{\TT}{Toro-Castro }
\newcommand{\HEOC}{Harten, Enquist, Osher, Chakravarthy }
\newcommand{\Network}{\mathcal{N}}
\newcommand{\Edges}{\mathcal{E}}
\newcommand{\Vertices}{\mathcal{V}}

\newcommand{\swh}{\ensuremath{h}} 	%height
\newcommand{\swH}{\ensuremath{H}} 	%absolute water level
\newcommand{\swq}{\ensuremath{q}} 	%massflux
\newcommand{\order}{k}
\newcommand{\x}{\ensuremath{x}}
\newcommand{\zeit}{\ensuremath{t}}
\newcommand{\edges}{\ensuremath{\mathcal{E}}}
\newcommand{\edge}[1]{\ensuremath{\edges_{#1}}}
	
\pgfkeys{/pgf/images/include external/.code={\includegraphics[draft=false]{#1}}}

\begin{frontmatter}
	
	
	
	%\title{ADER schemes and high order coupling on networks of hyperbolic conservation laws}
	\title{High order numerical methods for networks of hyperbolic conservation laws coupled with ODEs and lumped parameter models}
	%\author{Jochen Kall \corref{cor1}\fnref{fn1}}
	%\ead{kall@mathematik.uni-kl.de}
	%\address{Kurt-Schumacher-Stra\ss e 34, 67663 Kaiserslautern, Germany}
	%\fntext[fn1]{This is the first author footnote.}
	%
	%\author{Raul \corref{cor1}\fnref{fn1}}
	%\ead{kall@mathematik.uni-kl.de}
	%\address{Kurt-Schumacher-Stra\ss e 34, 67663 Kaiserslautern, Germany}
	%\fntext[fn1]{This is the first author footnote.}
	\author[JK]{\corref{cor1}Raul~Borsche}
	\ead{borsche@mathematik.uni-kl.de}
	\author[JK]{Jochen~Kall}
	\ead{kall@mathematik.uni-kl.de}
	\cortext[cor1]{Corresponding author}
	%\cortext[cor2]{Principal corresponding author}
	%\fntext[fn1]{This is the specimen author footnote.}
	%\fntext[fn2]{Another author footnote, but a little more longer.}
	%\fntext[fn3]{Yet another author footnote. Indeed, you can have
	%any number of author footnotes.}
	\address[JK]{Erwin Schr\"odinger Stra\ss e , TU Kaiserslautern, Building 48,
	67663 Kaiserslautern, Germany}
	
	\begin{abstract}
		In this paper we construct high order finite volume schemes on networks of hyperbolic conservation laws with coupling conditions involving ODEs.
		We consider two generalized Riemann solvers at the junction,
		one of \TT type and a solver of \HEOC type.
		The ODE is treated with a Taylor method or an explicit Runge-Kutta scheme, respectively. 
		Both resulting high order methods conserve quantities exactly if the conservation is part of the coupling conditions. 
		Furthermore we present a technique to incorporate lumped parameter models,
		which arise from simplifying parts of a network.		
		The high order convergence and the robust capturing of shocks is investigated numerically in several test cases. 
	\end{abstract}	 
	
	\begin{keyword}
		ADER 
		\sep
		Network
		\sep
		hyperbolic conservation law
		\sep
		WENO
		\sep
		generalized Riemann problem
		\sep
		Coupling
		\sep
		ODE
		\sep
		Runge-Kutta
		\sep
		Lumped Parameter Models
		
	\end{keyword}
	
\end{frontmatter}

%\linenumbers

\section{Introduction}
\label{S:Introduction}
Networks of hyperbolic PDEs arise from the modeling of many different problems, 
e.g. water and wastewater networks \cite{herty2012assessment,MR3281585,borsche2010coupling},
gas pipelines \cite{DomschkeKolbLang,brouwer_gas_2011,MR3335526},
traffic flow \cite{coclite2005traffic,MR3149318}, 
simulation of blood flow \cite{muller2014global,QuarteroniBlood,QuarteroniBlood2} 
or cell migration \cite{MZA:9150035}.
The description of such networks is based on one dimensional conservation laws along the edges 
and suitable coupling conditions at the nodes.
The simplest type of coupling uses a set of algebraic relations routing the flow between the arcs of the network.
In many of the above applications further coupling conditions arise in which an ODE is located in the junction
e.g. buffers \cite{HertyKlarPiccoli}, storage tanks,manholes \cite{DomschkeKolbLang,borsche2010coupling,MR3281585}
or the heart \cite{CNM:CNM2622}.
A  wide class of such coupling also occurs when so called lumped parameter models are applied to parts of the network \cite{CNM:CNM2622,QuarteroniBlood,QuarteroniBlood2,DomschkeKolbLang}. 
These models arise from simplification of the flow on the edges in regions where a coarser modeling can be afforded.

In all these applications fast and accurate numerical methods are needed.
For the one dimensional flow along the edges a huge variety of classical solvers for hyperbolic conservation laws is available \cite{toro2009riemann,leveque_finite_2002}.
Especially numerical methods of high order accuracy, e.g. WENO, ADER and DG schemes \cite{jiang1996efficient,toro2009riemann,hesthaven_nodal_2008},
achieve remarkably accurate solutions relative to their computational costs. 
For the application on networks however mainly first order schemes have been developed. 
Recently a high order Riemann solver for purely algebraic coupling conditions was presented in \cite{borsche2014ader}. 

In the present article we introduce two approaches to high order methods for vertices 
that can involve ODEs in addition to algebraic coupling conditions.
The first method follows the \TT approach of \cite{castro2008solvers,borsche2014ader}.
After solving a classical first order coupling problem
linear coupling conditions for the temporal derivatives of the states at the junction are considered.
The ODE is incorporated into this procedure 
by inserting the full Taylor-expansion of its solution into the coupling conditions.
The second method adapts the \HEOC approach \cite{castro2008solvers},
which solves a series of classical nonlinear Riemann problems.
These problems are considered at the supporting points in time of an explicit Runge Kutta scheme,
which is used for the discretization of the ODE.
In the context of this solver we investigate an efficient high order solver for lumped parameter models.
Here we aim to exploit the underlying network structure for the numerical method.

This paper is organized as follows:
First we formulate the problem, specify the coupling conditions % in \ref{S:Problem}
and define the generalized Riemann problem at a junction. % in \ref{S:Generalized Riemann Problem at a junction}. 
In section \ref{S:Generalized Riemann Problem at a junction} we recall the first order solver for such a coupling problem.
The high order method of \TT type is presented in Section \ref{S:GRP-TT},
the HEOC type solver in section \ref{S:GRP-HEOC}.
Based on these we describe a modification suited to networks including lumped parameter models.
Finally we present convergence studies and several numerical examples in section \ref{S:Numerik}
to show the quality of the schemes presented.

% We then develop the generalized Riemann solvers of \TT type in \ref{S:GRP-TT}  and of \HEOC type in chapter \ref{S:GRP-HEOC}.
% In chapter \ref{S:Lumped Parameter Models} we describe a GRP solver connecting to a lumped parameter network. 
% Finally we present convergence studies and some additional numerical demonstrations in \ref{S:Numerik}.

\section{Formulation of the problem and coupling conditions}
\label{S:Problem}
% We follow the same notation as we used in \cite{borsche2014ader} which we briefly repeat
% to consider a Network $\Network=(\Edges,\Vertices)$, consisting of edges connecting vertices 
A Network $\Network=(\Edges,\Vertices)$ consists of a set of edges $\Edges =\{E_1,\hdots,E_{\tilde n} \}$ which connect the vertices of the set $\Vertices=\{V_1,\hdots,V_{\tilde m} \}$. 
% \begin{align*}
% \Edges & =\{E_1,\hdots,E_{\tilde n} \}\ , &   
% \Vertices&=\{V_1,\hdots,V_{\tilde m} \}\ .
% \end{align*}
On each edge $E_i$, $i=1,\dots,\tilde n$ the quantities $u^i(x,t)\in \mathbb{R}^{d_i}$ are governed by a hyperbolic conservation law of the form
\begin{align}\label{eq:ConservationLaw}
\partial_t u^i + \partial_x f^i(u^i) & = 0\ 
\end{align}
with the flux function $f^i:\mathbb{R}^{d_i}\rightarrow \mathbb{R}^{d_i}$, time $t\in \R^{+}$ and location $x\in \left[0, L_i \right]$.

In the following we consider a single vertex $V$ and assume all edges to be oriented outwards, as depicted in  \figref{fig:Edge orientation convention}. 
% We assume all edges to be oriented outwards to formulate the coupling conditions, see \figref{fig:Edge orientation convention}.
\begin{figure}
	\centering
	
	\externaltikz{coupling-illustration}{
		\input{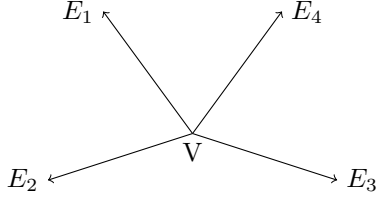}
	}		
	\caption{Edge orientation convention}
	\label{fig:Edge orientation convention.}
\end{figure}
Starting from this setup, networks of arbitrary shape can be easily constructed by elementary transformations, e.g. \cite{doi:10.1137/080716372}.
To improve readability we drop the index $j$ for all quantities at the junction, i.e. $w=w^j$, 
and the spatial dependencies of the $u^i$s, 
which are evaluated at $x=0$.

At the vertex $V$ we now assume a coupling of mixed algebraic-ODE type, i.e.
\begin{align}\label{eq:CouplingConditions}
% \Phi(u^1(t),\hdots,u^n(t),w)  =0\ , &&\text{(Algebraic Coupling Conditions)}\\
% \dot{w}=F(u^1(t),\hdots,u^n(t),w), &&\text{(ODE)}\\
	\begin{split}
	\Phi(u^1(t),\hdots,u^n(t),w)  =0\ , \\
\dot{w}=F(u^1(t),\hdots,u^n(t),w), \\
u^i(t)=u^i(0,t)\ .
	\end{split}
\end{align}
The algebraic coupling conditions are given by the function $\Phi:\bigotimes_{i=1}^n \mathbb{R}^{d_i} \times \R^{l} \rightarrow \mathbb{R}^{c}$ for $n$ connected edges 
and the ODE is defined by the flux $F:\bigotimes_{i=1}^n \mathbb{R}^{d_i} \times \R^{l} \rightarrow \R^{l}$.
Coupling conditions of that type arise from the modeling of e.g. storage components as manholes \cite{borsche2010coupling} or reservoirs \cite{MR3338429}, queues \cite{MR2368968,MR3050287}, 
as well as from representing parts of complicated networks by so called lumped parameter models \cite{QuarteroniBlood,QuarteroniBlood2,DomschkeKolbLang}.

In order to keep the number of coupling conditions constant,
the eigenvalues of the Jacobians $\nabla_{u^i} f^i$ $\lambda^i_l,\ l=1,\dots,d_{i}$ need to be bounded away from zero for all states considered
\begin{align}\label{eq:Eigenvalues}
\lambda^i_1 \leq \hdots \leq \lambda^i_{d_i}  \ ,
&   &   
|\lambda^i_l|>\tilde \epsilon>0 \qquad \forall l=1,\dots,d_{i}\ .
\end{align}
To ensure that the correct number of coupling conditions is provided and that the problem is well posed \cite{borsche2011mixed}, we require:
\begin{align}\label{eq:determinant}
\det\left( \nabla_{u^1}\Phi(u_g^1,\hdots,u_g^n,w_0) R^{+,1}|\dots|\nabla_{u^n}\Phi(u_g^1,\hdots,u_g^n,w_0) R^{+,n}\right) \neq 0\ , 
\end{align}
where $R^{+,i}=\left[r^i_{d_i-c_i+1}|\dots|r^i_{d_i}\right]$ denotes the matrix of all eigenvectors  of $\nabla_{u^i} f^i$ which belong to positive eigenvalues.
\begin{align}\label{eq:NumberoPositiveEigenvalues}
c_i =\#\{\lambda^i_j|\lambda^i_j>0\} 
\end{align}
% $c_i$
denotes the number of positive eigenvalues on edge $i$
and $c = \sum_{i=1}^n c_i$ defines the total number of coupling conditions prescribed by $\Phi$.

\section{The generalized Riemann problem at a junction}
\label{S:Generalized Riemann Problem at a junction}
One central building block for the construction of schemes in the ADER framework is the generalized Riemann problem.
In order to develop similar high order methods for networks, 
a detailed understanding of the generalized Riemann problem at the junction is required.
If additionally an ODE is located at the node, 
we aim to split the problem into two separate ones.
On the PDE side we are looking for high order approximations to the Godunov states at the boundaries of the PDE domains, while simultaneously evolving the ODE in the vertex one time step.

In the following we will discuss two variants to tackle this problem.
The first is based on the classical ADER approach of \TT accompanied by a Taylor method for the ODE.
The second one utilizes a \HEOC solver for the PDE and a Runge-Kutta scheme for the ODE.
% Constructing an ADER scheme boils down to solving generalized Riemann problems, a process we need to generalize even further to handle the coupling points with high order. 
% This splits in to separate problems. 
% Obtaining a high order approximation to the godunov states over time to construct high order fluxes for the PDEs at the boundaries and evolving the ODE in the vertex one time step. 
% We cover two variants to tackle this problem, 
% the first based on the classical ADER approach of \TT and a Taylor method for the ODE, 
% the second one based on the \HEOC solvers and a Runge-Kutta time stepping of the ODE.
\begin{Definition}
\label{def:generalized Riemann probem at a junction}
Generalized Riemann problem of order $k$ at a junction:\\
Consider a algebraic-ODE type coupling \eqref{eq:CouplingConditions} of $n$ edges governed by \eqref{eq:ConservationLaw}.
We call such a coupling situation with given initial state of the ODE $w_0$ 
and polynomial Riemann data $u^i(x,0)$ of order $k$ a Generalized Riemann problem of order $k$ at the junction

% We call a coupling situation with $n$ edges and coupling conditions of algebraic-ODE type 
% with given ODE state $w_0$ and polynomial Riemann data $u^i(x,0)$ of order $k$ a Generalized Riemann problem of order $k$ at the junction
\begin{align*}
\begin{split}
\Phi(u^1(t),\hdots,u^n(t),w(t))  &=0\ , \\
\dot{w}&=F(u^1(t),\hdots,u^n(t),w(t)),\\
u^i(t)&=u^i(0,t) \\
u^i(x,0)&=\sum\limits_{l=0}^{k-1} p^i_l\frac{x^l}{l!} \\
w(0)&=w_0\ .
\end{split}
\end{align*}
Analogously to the classical Riemann problem,
the states at the left boundary of the coupled edges $u_g^i(t)=\lim_{\tau\rightarrow 0+} u^i(0,t+\tau)$ are called Godunov states. 
\end{Definition}

\subsection{Solving the classical Riemann problem at the junction}
\label{S:Classical Riemann problem at the junction}
Before considering the generalized Riemann problem we investigate the classical Riemann problem at a junction with an ODE,
% An important building block we need later on is to solve the classical Riemann problem at an ODE junction 
i.e. the setup of definition \ref{def:generalized Riemann probem at a junction} with $k=1$
\begin{align*}
\Phi(u^1(t),\hdots,u^n(t),w(t))  &=0\ , &
\dot{w}&=F(u^1(t),\hdots,u^n(t),w(t)),\\
u^i(t)&=u^i(0,t) &
u^i(x,0)&=p^i_0\\
w(0)&=w_0\ .
\end{align*}
In the case without ODE, the Godunov states are constant in time, i.e. $u^i_g(t)\equiv u^i_g(0)$ $\forall t\geq0$.
Due to the presence of the ODE the state in the junction can vary over time and thus also the Godunov states $u^i_g(t)$ change.
However for solving the generalized Riemann problem we are only interested in the states at $t\rightarrow 0+$.

Since the solution of the ODE is continuous in time \cite{borsche2011mixed} we have $\lim_{t\rightarrow 0+}w(t)=w_0$.
Knowing the initial state of the ODE, the problem at $t=0+$ reduces to a classical Riemann problem at a junction.
This we can solve with the help of the so called Lax-Curves \cite{borsche2014ader,colombo2008cauchy}.

Solving such a classical Riemann problem at a junction is equivalent to finding a set of states that fulfill the algebraic part of the coupling conditions
while being possible Godunov states accessible from the Riemann data on each edge.
In order to be such an accessible state they have to lie on the concatenated Lax curves anchored in the right states $u^i_r=p_0^i$ \cite{colombo2008cauchy}:
\begin{align*}
L_g^i(\xi^i_1,\hdots,\xi^i_{c_i},u^i_r)
& =L_{d_i-c_i+1}(\xi^i_1,\cdot)\circ \hdots \circ L^i_{d_i}(\xi^i_{c_i},u^i_r)\ , 
\end{align*}
where the number of curves and free parameters $c_i$ is given by \eqref{eq:NumberoPositiveEigenvalues}.
% $c_i=\#\{\lambda^i_j|\lambda^i_j>0\}$. 
The operator $\circ$ denotes the concatenation in the last variable, i.e for two functions $g$ and $h$
\begin{align*}
g(\xi_1,\hdots,\xi_l,\cdot)\circ h(\xi_{l+1},\hdots,\xi_m,x) & =g(\xi_1,\hdots,\xi_{l},h(\xi_{l+1},\hdots,\xi_m,x))\ . 
\end{align*}
Therefore we have to solve the equations
\begin{align}\label{eq:0thOrderCoupling}
\Phi
\left(
L^1_g(\xi^1,u^1_r),L^2_g(\xi^2,u^2_r),\hdots,L^n_g(\xi^n,u^n_r),w_0	 	
\right)
& =0 
\end{align}
for the unknowns $\xi^i=(\xi^i_1,\hdots,\xi^i_{c_i})$.
The local solvability of this system for states close to $u_r$ is assured by condition \eqref{eq:determinant} \cite{colombo2008cauchy}.
Once the parameters $\xi^i$ are known, the Godunov states can be determined by evaluating the concatenated Lax-curves
\[ u^i_g(0)=L_g^i(\xi^i,u_r^i)\ .\]
Thus the complete set of states at the junction at $t= 0+$ is given by $\{u^1_g(0),\hdots,u^n_g(0),w_0\}$.
% Furthermore, with the Godunov states known, $\dot{w}(0)=F(u^1_g(0),\hdots,u^n_g(0),w(0))$ can be evaluated, a quantity we need in our schemes.

\section{Generalized Riemann solver of \TT type}
\label{S:GRP-TT}
In the ADER framework the \TT approach provides a procedure to construct a polynomial in time approximating the solution of the generalized Riemann problem at the considered interface.
This is achieved by splitting up the problem into one classical non linear Riemann problem 
and $k-1$ linearized Riemann problems for the temporal derivatives of the Godunov states. 
An extension of this procedure to junctions without ODEs was presented in \cite{borsche2014ader}.

Let a generalized Riemann problem at an ODE junction be given as in definition \ref{def:generalized Riemann probem at a junction}.
As a first step we solve the zeroth order classical Riemann problem as described in section \ref{S:Classical Riemann problem at the junction}
using the zero order data, i.e. $u_r^i=p^i_0$ .
Note that once the states at $t=0+$ are known, we can directly evaluate the ODE 
$\dot{w}(0)=F(u^1_g(0),\hdots,u^n_g(0),w(0))$, 
which already provides some information about the development of $w(t)$.
% , a quantity we need in the next step and to propagate the ODE.
% This construction guarantees that the states at the junction are reached only by waves traveling into the edges and no information is lost due to waves exiting the domain.

In a second step we aim to compute the temporal derivatives of the involved states.
As in the classical ADER framework,
we obtain governing equations for the derivatives by differentiating the conservation laws with respect to $t$ and obtain
% 
% To obtain Riemann problems for the derivatives of the states, we differentiate the conservation laws with respect to $t$ and get:
\begin{align}\label{eq:generalizedRP}
\partial_t(\partial_t^k{u^i}) + \nabla f^i(u_g^i)\partial_x(\partial_t^ku^i) +\text{'sources'} =0 &   & k = 1,\dots,{k_{max}}\ . 
\end{align}
The term 'sources' encompasses everything that only depends on derivatives of degree $k$ and less.
It can be dropped, since it does not immediately act on the Godunov states \cite{toro2009riemann}.
These governing equations for the temporal derivatives are linear hyperbolic systems. 
Therefore the corresponding Lax curves  are linear as well.
The concatenated Lax curves to a Riemann problem  with $\partial_t^k u_r^i$ as states on the right hand side
have the short form
\begin{align}\label{eq:linearLaxCurves}
L^{i,k}_g(\xi^i,u_r)=\partial_t^k u_r^i+R^{+,i}\xi_k^i\ ,
&   &   
R^{+,i}=\left[  
r^i_{d_i-c_i+1} | \hdots | r^i_{d_i}
\right]
\ ,
\end{align}
where $r^i_j$ denotes the eigenvector corresponding to the $j$-th eigenvalue $\lambda_j^i$ of $\nabla f^i(u^i_g)$.
Note that \eqref{eq:Eigenvalues} still holds as the Jacobian is the same as in the Riemann problem of order zero.

In order to obtain the coupling conditions for the temporal derivatives, we differentiate $\Phi$ with respect to time.
The first order derivative of $\Phi$ reads
\begin{align*}
% \Phi(u^1,\hdots,u^n,w)                       & =0     \\
% \Rightarrow \partial_t\Phi(u^1,\hdots,u^n,w) & =0     \\
% \Rightarrow 
\frac{d}{dt}\Phi(u_g^1,\hdots,u_g^n,w) & =0     \\
\Rightarrow \sum_{l=1}^n \nabla_{u_g^l}\Phi(u_g^1,\hdots,u_g^n,w)
\partial_tu_g^l                
+\nabla_w  \Phi(u_g^1,\hdots,u_g^n,w)\partial_tw            & =0 \ , 
\end{align*}
where all quantities are evaluated at $t=0$ and $u_g^i$ denotes the Godunov state at $x=0$ on edge $i$.
By inserting the ODE $\partial_tw=F(u_g^1,\hdots,u_g^n,w)$ we obtain 
\begin{align}\nonumber
\nabla_{u}\Phi(u_g^1,\hdots,u_g^n,w)
\partial_t u_g                
+\nabla_w  \Phi(u_g^1,\hdots,u_g^n,w)F(u_g^1,\hdots,u_g^n,w)            & =0 \\ 
\label{eq:first_order_derivative_CC}
\nabla_{u}\Phi(u_g^1,\hdots,u_g^n,w)
\partial_tu                
+\Psi_1(u_g^1,\hdots,u_g^n,w)            & =0 \ .
\end{align}
The function $\Psi_1$ only depends on states at $t=0$ and not on any temporal derivative. 
Thus we have a linear system governing the temporal derivatives of the Godunov states at the junction. 
Here and in the following we assume the coupling conditions $\Phi$ to be sufficiently differentiable, 
otherwise we consider the usage of high order schemes not appropriate.

For the derivatives of orders $k\geq 2$ additional terms arise, which we summarize as $\tilde\Psi_k$
\begin{align*}
\nabla_u\Phi(u_g^1,\hdots,u_g^n,w)\partial_t^ku_g
+\tilde\Psi_k(u_g,\partial_tu_g,\hdots,\partial_t^{k-1}u_g,w,\partial_tw,\hdots,\partial_t^{k}w)=0\ .
\end{align*}
Note that these lower order terms can not be dropped as in the classical ADER framework, i.e. equation \eqref{eq:generalizedRP},
since we can not expect that it acts delayed in any form. 
It is easy to see that $\tilde\Psi_k$ only depends on lower order derivatives of $u_g$
but still on all $k$ derivatives of $w$.
This we can reduce by using the derivative of the ODE itself
\begin{align}\label{eq:derivativeODE}
\begin{split}
\partial_t^{k}{w}&=\partial_t^{k-1}F(u^1_g(t),\hdots,u^n_g(t),w) \\
&=\Xi(u^1_g,\partial_tu^1_g,\hdots,\partial_t^{k-1}u^1_g,\hdots,u^n_g,\partial_tu^n_g,\hdots,\partial_t^{k-1}u^n_g,
w,\partial_tw,\hdots,\partial_t^{k-1}w)\ .
\end{split} 
\end{align}
% In fact any derivative of $w$ can be expressed in this manner.
Thus we end up with a linear equation for $\partial_t^ku_g$ which only requires derivatives of lower order than $k$
\begin{align}\label{eq:DiffCouplingConditions}
\nabla_u\Phi(u_g^1,\hdots,u_g^n,w)\partial_t^ku_g 
+\Psi_k(u_g,\partial_tu_g,\hdots,\partial_t^{k-1}u_g,w,\partial_tw,\hdots,\partial_t^{k-1}w)=0\ .
\end{align}
Finally we can use this expression to obtain all the temporal derivatives of $u_g$ iteratively,
by starting with the first order equation \eqref{eq:first_order_derivative_CC} 
and successively solving \eqref{eq:DiffCouplingConditions} in increasing order of $k$.

To solve each of these systems \eqref{eq:DiffCouplingConditions},
we proceed as in the zeroth order case.
The temporal derivatives of $u_g$ are governed by a linear conservation law \eqref{eq:generalizedRP}
and coupled by a set of linear coupling conditions \eqref{eq:DiffCouplingConditions}.
Thus we insert the concatenated linear Lax-curves \eqref{eq:linearLaxCurves} with the free parameter $\xi^i_k$
% Plugging in the lax curves yields:
\begin{align*}
\sum_{i=1}^n \nabla_{u^i}\Phi(u_g^1,\hdots,u_g^n,w)	
\left(
R^{+,i}\xi_k^i+\partial_t^ku^i_r 
\right) +\Psi_k           & =0     
\\
\Rightarrow \sum_{i=1}^n \nabla_{u^i}\Phi(u_g^1,\hdots,u_g^n,w)	
R^{+,i}\xi_k^i 
+\sum_{i=1}^n \nabla_{u^i}\Phi(u_g^1,\hdots,u_g^n,w)
\partial_t^ku^i_r +\Psi_k & =0 \ . 
\end{align*}
By introducing the notations
\begin{align*}
a_i    & =\nabla_{u^i}\Phi(u_g^1,\hdots,u_g^n,w) R^{+,i} \ ,
&   &   &   
A      & =\left(                          
a_1 | a_2 | \hdots | a_n 
\right)\ , \\
\xi_k    & =                                
\begin{pmatrix}                 
\xi_k^1  &                                  
\xi_k^2  &                                  
\hdots &                                  
\xi_k^n 
\end{pmatrix}^T \ ,
&   &   &   
\partial_t^k u_r&=
\begin{pmatrix}
\partial_t^ku_r^1  &                                  
\partial_t^ku_r^2  &                                  
\hdots &                                  
\partial_t^ku_r^n 			
\end{pmatrix}^T\ ,
\end{align*}
this linear system can be written as
\begin{align}\label{eq:CouplingDerivatives}
% 	\Rightarrow 
A\xi_k+\nabla_u\Phi(u_g^1,\hdots,u_g^n,w)\partial_t^ku_r+\Psi_k & =0 \ . 
\end{align}
The matrix $A$ is exactly the one in \eqref{eq:determinant} considered for the well-posedness of the coupling conditions.
Since we have $\det A \neq 0$, we can solve for the unknowns $\xi_k$
\begin{align*}
% 	\Rightarrow 
\xi_k
& =A^{-1}\left( 
-\nabla_{u}\Phi(u_g^1,\hdots,u_g^n,w)\partial_t^ku_r
-\Psi_k
\right)\ .
\end{align*}
In order to evaluate this expression, 
the temporal derivatives of the states within the edges $\partial_t^ku_r$ are needed.
Analogous to the classical ADER approach we start with a WENO reconstruction of the spatial initial data.
Since we are at the boundary of a domain, we have to use an one-sided reconstruction of type \cite{tan2010inverse} to obtain spatial derivatives $\partial_x^ku_r$ at the junction. 
These we can transform into temporal derivatives $\partial_t^ku_r$ using the Cauchy-Kowalewski or Lax-Wendroff procedure \cite{toro2009riemann,tan2010inverse}.
%Note that it is important that the Cauchy-Kowalewski procedure is carried out on the basis of the zeroth order Godunov state,
%which we can obtain from the classical zeroth order Riemann Problem at the junction.

Note that carrying out the Cauchy-Kowalewski procedure on the basis of the zeroth order Godunov state is not strictly necessary, but does not cause any additional computational costs either. It is however necessary to make the scheme revert to the classical Toro-Titarev ADER scheme for one on one coupling, as we showed in \cite{borsche2014ader}.
Once the values of $\xi_k$ are determined, 
we can evaluate the Lax-curves \eqref{eq:linearLaxCurves} to obtain the temporal derivatives of the Godunov states $\partial_t^ku_g$.
With these we can directly build a polynomial approximation of $u_g(t)$ at the interfaces of the junction.

% The state of the ODE can be approximated directly.
An approximation to the states of the ODE can now be constructed easily.
Since for the computation of $\partial_t^ku_g$ we already determined all the derivatives of $w$ using equation \eqref{eq:derivativeODE},
we just combine these to a Taylor series approximating $w(t)$.

% Since all states on the right hand side are known from previous steps, we can directly compute $\partial_t^{k+1}{w}$.
% With all the $k_max$ derivatives of $w$ we can set up a Taylor series approximating the state of the ODE.
% Plugging into the linearized lax curves gives a Taylor series of the godunov states, 
% which in turn gives a $k-th$ order numerical flux, analogously, 
% the taylor series of the ODE state gives a $k-th$ order Taylor scheme. 
% \subsection{Filling Ghostcells and obtaining the Riemann data}
% This still works exactly as described in \cite{borsche2014ader}, we therefore skip it here.
% \subsection{Description of the numerical scheme}
% \label{S:scheme TT}
The scheme derived from the \TT type solver for the generalized Riemann problem at the junction with ODE can be summarized as follows:
\begin{enumerate}
%	\item Fill the ghost cells at the junction according to section \ref{sec:spatialreconstruction_junction}.
	\item Obtain GRP data at the junction via one-sided polynomial reconstruction.
	\item Solve the zeroth order Riemann problem at the junction, as described in section \ref{S:Classical Riemann problem at the junction}.%, to obtain the Godunov states and $\dot{w}$.
	\item Apply the Cauchy-Kowalewski procedure to obtain temporal polynomials as input data.
% 	\item Solve the generalized Riemann problem at the junction, as described in section \ref{S:GRP-TT}.	 
	\item Solve the generalized Riemann problem at the junction, as described above.	 
% 	\item Compute the fluxes across cell interfaces at the junction using the temporal derivatives of the Godunov states and the higher order derivatives of $w$.
	\item Approximate the fluxes across cell interfaces at the junction using the temporal derivatives of the Godunov states.
	\item Update the state $w$ in the junction by applying a Taylor scheme of order $k$.
	\item Fill the ghost cells at the junction if needed.
	\item Run a high order finite volume scheme e.g. ADER to compute the fluxes across interior cell interfaces.
\end{enumerate}
A detailed description how to fill the ghost cells, if needed, is given in \cite{borsche2014ader}.
% Possible ones-sided reconstructions are discussed in \cite{}.

At this point we note that the above procedure requires,
besides the \CK procedure,
the derivatives up to order $k$ of the ODE 
as well as those of the coupling conditions.

% Note that this solver has some advantages and disadvantages. 
% Conserved properties are conserved exactly if the conservation is part of the coupling conditions, 
% but the method requires the evaluation of the $k$-th derivative of ODE as well as coupling conditions, 
% which might be very inconvenient or downright impossible. 
% In This case we recommend the use of the \HEOC type solver described in \ref{S:GRP-HEOC}, 
% which does not have these requirements, but conserves quantities only up to order $k$.

\section{Generalized Riemann solver at an ODE junction \HEOC type}
\label{S:GRP-HEOC}
A popular alternative to the \TT approach for solving generalized Riemann problems is the \HEOC solver \cite{castro2008solvers}.
Instead of solving the linear governing equations for the derivatives,
it uses several classical Riemann problems at different points in time to achieve a high order approximation of the flux.
Here we present a HEOC type solver for junctions, which is accompanied by a Runge-Kutta scheme for the ODE.

\begin{figure}
\centering
\externaltikz{HEOC-schematic}
{
\input{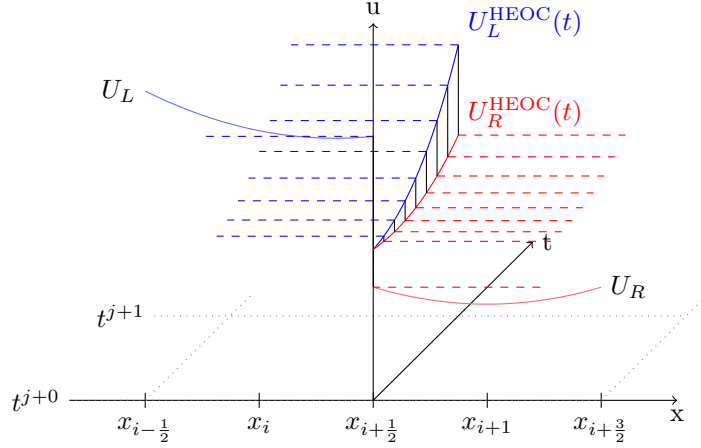}
}
\caption{HEOC schematic.}
\label{fig:HEOC schematic}
\end{figure}

% \subsection{HEOC approach}
In the classical HEOC approach, the solution of a generalized Riemann problem is approximated by the solutions of a sequence of classical Riemann problems. 
% These are located on the time axis at the supporting points of a given quadrature rule,
% which is used for the integration of the fluxes at the interface.
First we fix a quadrature rule according to the desired order of the scheme.
Then the classical Riemann problem at $t=0$ is solved.
Based on its solution the spatial data on each side of the interface can be flipped into the time domain using the \CK procedure.
These polynomials can be evaluated at the supporting points of the given quadrature rule,
such that a series of classical Riemann problems arises.
Their solutions serve as approximation of the states at the interface at these points in time.
This information can be inserted into the quadrature rule to obtain a high order approximation of the fluxes at the interface.
In \figref{fig:HEOC schematic} a schematic representation of this procedure is shown.

In the following we adapt this approach to generalized Riemann problems at junctions with ODEs. 
Before considering the fluxes of the PDE, we start with the numerical method for the ODE.
Here we choose an explicit Runge-Kutta scheme which is at least accurate of order $k$.
Usually the coefficients of RK-schemes are given in form of a Butcher array
% \begin{Definition}
% 	Explicit Runge Kutta scheme for an ODE
{
\renewcommand{\arraystretch}{1.2}
\begin{align*}
% \dot{y}&=F(t,y) &
\bar{B}&=
\left(
\begin{array}{c|c}
\bar{c} & \bar{A} \\ 
\hline
& \bar{b} 
\end{array}
\right) \ ,
\end{align*}
}
such that for an ODE $\dot{w}=F(t,w)$ the update formula reads
\begin{align}\label{eq:RKscheme_update}
w_{n+1}&=y_n+\Delta t\sum_{j=1}^{s}\bar{b}_jk_j \ ,&
k_i&=F(t_n+\bar{c}_i\Delta t,w_n+\Delta t\sum_{j=1}^{i-1}\bar{a}_{ij}k_j)\ .
\end{align}
% \end{Definition}
At this point we note that each RK-scheme of order $k$ naturally provides a quadrature formula of order $k$
with the supporting points $c_i\Delta t$, such that
\begin{align*}
 \int_0^{\Delta t}\tilde F (\tau)d\tau = \Delta t \sum_{j=1}^{s}\bar{b}_j\tilde F(t_n+\bar{c}_i\Delta t) +\mathcal{O}(\Delta t^{k+1})\ .
\end{align*}
In the following we will use exactly these intermediate time levels $t_l=c_l\Delta t$ to set up the HEOC coupling procedure.

Before considering higher order terms we have to solve the classical Riemann problem at $t=0$,
\begin{align*}
\Phi(p^1_0,\hdots,p^n_0,w_0)=0\ .
\end{align*}
Following the procedure described in \ref{S:Classical Riemann problem at the junction} we obtain the values of Godunov states at $t=0$, i.e.  $u^1_g(0),\hdots u^n_g(0)$. % and the first stage of the RK scheme.
Using these values we can transform the spatial data from a one-sided polynomial reconstruction into temporal data via \CK procedure.
Thus for each connected edge we obtain a temporal polynomial of the form
\begin{align*}
% 	u^{i,\text{HEOC}}(t)&=\sum\limits_{l=0}^{k-1} p^{i,\text{HEOC}}_l\frac{t^l}{l!} \ .
	u^{i}_r(t)&=\sum\limits_{l=0}^{k_{\max}-1} p^{i,\text{time}}_l\frac{t^l}{l!} \ 
\end{align*}
as input data for the generalized Riemann problem at the junction.
With these values available, we now aim to solve the classical Riemann problems at the time levels $t_l$
\begin{align}\label{eq:HEOC_Phi_wtl}
% 	\Phi(u^{1,\text{HEOC}}(t_l),\hdots,u^{n,\text{HEOC}}(t_l),w_l) &=0
	\Phi(u^{1}_g(t_l),\hdots,u^{n}_g(t_l),w(t_l)) &=0\ .
\end{align}
In order to apply the technique of section \ref{S:Classical Riemann problem at the junction} 
we have to provide some approximation for the value of $w(t_l)$.
% Since we can solve these equations only for the $u^i_g$, we have to provide some approximation to the value $w(t_l)$.
This is naturally provided by the RK-scheme as in the second formula of \eqref{eq:RKscheme_update}
\begin{align}\label{eq:w_l}
%      t_l&=c_l\Delta t &
	 w_l&=w_0+\Delta t \sum\limits_{i=1}^{l-1}\bar{a}_{l,i}k_i\ .
\end{align}
As in the classical RK methods, the value $w_l$ is not necessarily an approximation of very high order, 
but chosen in such a way that in the final update the desired order is obtained. 
Note that for the evaluation of the stages $k_i=F(t_i,u_g^1 (t_i),\hdots, u_g^n(t_i),w_i)$ $ i=1,\dots,l-1$ values of $w_i$ and $u_g(t_i)$ are needed. 
But since we have chosen an explicit RK-scheme, only data from the previous $l-1$ stages is used.
% All values which are needed to evaluate the stages $k_i=F(u_g^1 (t_i),\hdots, u_g^n(t_i),w_i)$ $ i=1,\dots,l-1$, which depend on $w$ and $u$, are known from the previous classical Riemann problems.

Inserting \eqref{eq:w_l} into \eqref{eq:HEOC_Phi_wtl} we obtain
\begin{align*}
% 	\Phi(u^{1,\text{HEOC}}(t_l),\hdots,u^{n,\text{HEOC}}(t_l),w_l) &=0
	\Phi(u^{1}_g(t_l),\hdots,u^{n}_g(t_l),w_l) &=0\ 
\end{align*}
and can solve this classical Riemann problem for the Godunov states at $t_l$.

Once all Riemann problems are solved successively, we can compose the solutions to determine the fluxes of the conservation laws and to update the ODE
\begin{align}\label{eq:ODEupdate_fluxintegration}
f^i_{-\frac{1}{2}}&=\sum\limits_{l=1}^{k} \bar{b}_l f^i(u^i_g(t_l))\ ,&
w(t+\Delta t)&=w(t)+\Delta t \sum\limits_{l=1}^{k}\bar{b}_l k_l\ .
\end{align}
The coefficients $b_l$ in both formulas are those of the RK-scheme \eqref{eq:RKscheme_update}
and the stages $k_l$ have been already computed for the formula \eqref{eq:w_l}.
% Finally we use the Quadrature rule of the RK scheme to approximate the value of $w$ at $t+ \Delta t$ and the fluxes over the boundaries of the edges:
The complete scheme can be summarized as
\begin{enumerate}
	\item Obtain GRP data at the junction via one-sided polynomial reconstruction.
	\item Solve the zeroth order Riemann problem at the junction, as described in section \ref{S:Classical Riemann problem at the junction}.%, to obtain the Godunov states and $\dot{w}$.
	\item Apply the Cauchy-Kowalewski procedure to obtain temporal polynomials as input data.
	\item Solve the classical Riemann problems at the times $t_l$ as described above.
	\item Approximate the fluxes across cell interfaces at the junction using the Godunov states at the time levels $t_l$.
	\item Update the state $w$ in the junction by applying the RK-scheme of order $k$.
	\item Run a high order finite volume scheme to compute the fluxes across interior cell interfaces.
\end{enumerate}

One advantage of this approach is that neither the derivatives of the coupling conditions nor those of the ODE are needed.
Thus the only symbolic manipulation necessary is the CK procedure, which is required for ADERs scheme inside the domain anyway. 
Furthermore we can solve the ODE with some classical RK-scheme, 
which is helpful especially for complicated or large ODEs
e.g. those that arise from lumped parameter models.

The main disadvantage of this approach is the higher computational costs, 
since several nonlinear Riemann problems have to be solved instead of just one nonlinear and a couple of linear ones.
Furthermore we do not have enough information at the interfaces to fill possible ghost cells at the computational boundary. 
Thus for the reconstruction in the interior of the domain but close to the boundary we use a reconstruction with variable stencil lengths \cite{tan2010inverse}.

\subsection{Conservation of quantities}
In many applications and their corresponding models some quantities are conserved in the complete network, e.g. the total mass \cite{MR3050287,MR3281585,MR3335526,QuarteroniBlood}.
This is not only established via the conservation laws on the edges, 
but also due to a careful choice of the coupling conditions and the ODE in the junction.

Since the conservation is guaranteed to be exact in the interior of the edges by the numerical scheme, 
it is desirable that also the coupling procedure is conservative.
The first order junction solver in section \ref{S:Classical Riemann problem at the junction} is conservative,
since the same Godunov states are used in the coupling conditions 
as well as for the computation of the fluxes across the interfaces.
If the ODE is updated by an explicit Euler scheme, 
its update relies on the same values as the coupling condition.
Clearly at the junction the conservation is only guaranteed up to the precision of the numerical method used to solve the nonlinear system arising from the coupling conditions.

In case of no ODE in the junction it has been proven in \cite{borsche2014ader} that the \TT approach also conserves the selected quantities.
This proof can be easily modified such that it fits the current setting.
We just have to take care that the ODE is updated with exactly the same numerical values as those arising in the coupling procedure.
Therefore it is mandatory to use a Taylor scheme for the ODE.

% Here we shortly motivate how quantities, which are conserved in the coupled system of PDEs and ODE, are also conserved by the proposed numerical methods.
% For the \TT approach the prove of \cite{borsche2014ader} can be easily modified to fit to the current setting.
For the HEOC approach at any intermediate time level $t_l$ a classical Riemann Problem at a junction is considered.
If now some quantity is conserved in the underlying system, at each of these Riemann problems the fluxes of the resulting Godunov states and the flux of the ODE balance exactly. 
Since we choose the identical $\bar{b}_l$s for the flux integration and the RK-scheme \eqref{eq:ODEupdate_fluxintegration}, this also holds for the final updates in the PDEs and the ODE. 

% The big advantage of this approach is that no symbolic manipulation beyond the CK procedure, 
% which is already required for ADER schemes anyways, is necessary. 
% This especially comes in handy when dealing with ODE vertices with an variable number of components, 
% e.g. lumped parameter models.
% \todo{conservation !!!}

\subsection{Source-terms}
In many applications the conservation law \eqref{eq:ConservationLaw} is replaced by a balance law via introducing source terms.
These usually do not affect the coupling procedure \cite{MR3281585}.
If we have a numerical scheme at hand that is capable to treat the source terms properly 
and include the sources into the \CK-procedure, all the above methods can be applied.
Since the lower order terms in \eqref{eq:generalizedRP} are dropped, 
the sources do neither change the Lax curves nor the governing equations of the higher order derivatives.

\section{Lumped Parameter Models}
\label{S:Lumped Parameter Models}
In a lot of real world applications the dimensions of the network exceed the affordable computational effort, e.g. capillaries in the circulatory system.
At the same time a detailed description of the flow is only needed in certain areas of the network.
Therefore it is often convenient to describe some parts of the network by simpler models.
A wide class of such reduced models are the so called lumped parameter models,
which are used to describe e.g. the human circulatory system \cite{QuarteroniBlood,QuarteroniBlood2} or gas networks \cite{DomschkeKolbLang}.

% A typical situation where a network of PDEs is coupled with ODEs at a node or the boundary arises when so called lumped parameter models are used.
% These result from a simplification of parts of the network which do not require detailed modeling or are too big to simulate entirely. 
% Lumped parameter models play an important role in the applications of hyperbolic networks \cite{QuarteroniBlood,QuarteroniBlood2}.
In this section we explain a process to construct high order schemes for hybrid models of networks
containing hyperbolic conservation laws and lumped parameter models.
To have access to such a process in an algorithmic framework is especially of interest in the context of dynamical switching between highly resolved and reduced models \cite{DomschkeKolbLang}.

% We then show the designed high order of convergence in \ref{S:Numerik}.
% 
Consider one edge in the network of length $L$ with a conservation law
\[\partial_tu(t,x)+\partial_x\left(f(u(t,x))\right)=0\ .\] 
Following the approach proposed in \cite{QuarteroniBlood2}, a lumped parameter model is obtained by averaging over the whole spatial domain, i.e.
\begin{align}
% \frac{1}{L}
% \int_{0}^{L}
% \left(
% \partial_tu(t,x)+\partial_x\left(f(u(t,x))\right)
% \right)d_x &=0  \\
% \Rightarrow
\label{eq:LumpingPDE}
\partial_t
\left(
\underbrace{
\frac{1}{L}
\int_0^L
u(t,x) d_x
}_{=:U}
\right)
&+\frac{1}{L}
\left[
f(\underbrace{u(t,L)}_{=:U_g^r})-f(\underbrace{u(t,0)}_{=:U_g^l})
\right]=0
\\
\Rightarrow
U'&=F(U,U_g^l,U_g^r)\ . \nonumber
\end{align}
Thus the averaged state $U$ is governed by a simple ODE with the unknown Godunov states at the left and right boundary $U_g^l$ and $U_g^r$.
If such an edge is connected to a node in the network, these values will be determined when solving the associated coupling conditions. 
The coupling procedure from section \ref{S:Classical Riemann problem at the junction} remains unchanged by the averaging process,
but the Lax curve have now to be anchored at the only accessible state $U$. 
This implies, if complete sections of the network are simplified to lumped parameter models in the above manner, 
that the PDEs are not coupled to a simple ODE but to a differential algebraic equation (DAE).
% i.e. an Differential Algebraic Equation (DAE) coupled with the adjacent vertices in the usual way.
% Performing this procedure for a sub-network results in a bigger set of DAEs, i.e.\\
The components of the lumped region can be summarized as follows
\begin{align}
\label{eq:lumped_ODE}
w' &= F(w,u_g) &\text{\parbox{8cm}{ODE part, originates from lumped edges and ODE parts already present in vertices,}} \\
\label{eq:lumped_coupling_conditions}
\Phi(u_g,w)&=0 &\text{\parbox{8cm}{Algebraic constraints stem from the coupling conditions of the vertices,}} \\
% u_g-L(w,u_R,\xi)&=0 & \text{\parbox{8cm}{Lax curve condition to restrict the solutions to valid godunov states}}
u_g-L(w,\xi)&=0 & \text{\parbox{8cm}{Lax curve condition to connect the Godunov states with the internal states.}} \nonumber
\end{align}
Note that here $w$ contains the averaged states of all lumped edges in a connected area 
and all possible ODE states in the junctions of this region.

From a numerical point of view the algebraic constraints \eqref{eq:lumped_coupling_conditions}
do enforce the usage of an appropriate solver of DAEs in the junction, e.g. modified RK schemes \cite{MR2657217}.
% There are modified RK schemes to deal with DAEs in general, i.e. \todo{Referenzen, namen Heirer oder?},
In the following however we  want to present a technique to construct a high order solver using the underlying network structure of the LPM.
% \todo{Vorteile unser Ansatz}

Since the lumped areas might be very large or vary in time, due to switching between the models,
we base the following procedure on the HEOC approach.
% In the first step we solve the coupling conditions in all junctions of the network,
% which provides sufficient data to proceede to the next stage of the underlying RK-scheme.
% As for the next stage of the ODE \eqref{eq:lumped_ODE} the values $u_g$ in the junctions are needed,
% but canbe obtained by solving again the coupling conditions as described in section \ref{S:GRP-HEOC}.
% % The flux function of the ODE obviously depends on the godunov states at the ends of the edges, which still are subject to the coupling conditions of the adjacent vertices. 
% % This means that solving the coupling conditions of a LPM model of this type can be done on a vertex by vertex basis exactly as it is done for normal PDE edges. 
% % This allows us to adapt the HEOC approach of \ref{S:GRP-HEOC}
% % by using the stages of of the RK scheme for interior edges and ODE vertices as intermediary values for the next Runge Kutta stage.
As in section \ref{S:GRP-HEOC} we first apply polynomial reconstruction to obtain spatial data on the PDE-edges.
For each vertex in the lumped parts of the network we solve the following zeroth order coupling problem as \eqref{eq:0thOrderCoupling}
\begin{align*}
\Phi(u^1_g,\hdots,u^i_g,w_0) =0&
&\text{ with }
&&u^i_r= \begin{cases}
 	p^i_0  & \text{edge is equipped with a PDE} \\
 	U^i_0  & \text{edge model is lumped}
\end{cases} \quad.
\end{align*}
This yields the Godunov states on the PDE edges at $t=0$, which are used to flip the data into time by \CK procedure,
providing polynomials $u^{i,\text{time}}(t)$.
We repeat this step for each stage $l$ of the RK-scheme in the HEOC approach and for each vertex in the lumped network,
i.e. we solve the zeroth order coupling problem at $t_l=c_l\Delta t$,
\begin{align*}
\Phi(u^1_g,\hdots,u^i_g,w_l) =0&& \text{ with }&&
u^i_r&= \begin{cases}
u^{i,\text{time}}(t_l)  & \text{If the edge in question is a PDE} \\
U^i_l  & \text{if the edge in question is lumped}
\end{cases} \quad .
\end{align*}
The states $U^i_l$ are known from the previous stage of the RK-scheme
and the new ones are obtained by
\begin{align*}
	U^i_{l+1}=U_0+\Delta t \sum\limits_{i}^{l}\bar{a}_{l+1,i}k^U_i\ ,
&&
w^j_{l+1}=w_0+\Delta t \sum\limits_{i}^{l}\bar{a}_{l+1,i}k^w_i\ ,
\end{align*}
where $k_i^U$ and $k_i^w$ are the components of the intermediate values $k_i$ for the lumped edges or the ODEs in the junctions respectively.

With this procedure we can keep the full network structure, 
such that we can easily select any part of this network and switch between the simplified and the more accurate model.
The scheme can be summarized by the following steps
\begin{enumerate}
	\item Obtain GRP data at the junction via one-sided polynomial reconstruction.
	\item Solve the zeroth order Riemann problem at the junction, as described in section \ref{S:Classical Riemann problem at the junction}.%, to obtain the Godunov states and $\dot{w}$.
	\item Apply the Cauchy-Kowalewski procedure to obtain temporal polynomials as input data.
	\item Compute each stage of the RK-scheme as described above.
	\item Approximate the fluxes across cell interfaces at the junction using the Godunov states at the time levels $t_l$.
	\item Update the state $w$ in the junction by applying the RK-scheme of order $k$.
	\item Run a high order finite volume scheme to compute the fluxes across interior cell interfaces.
\end{enumerate}

\subsection{Source-terms in lumped parameter models}
% In some applications the conservation law \eqref{eq:ConservationLaw} is replaced by a balance law via introducing some source terms.
% These usually do not affect the coupling procedure \cite{MR3281585}.
% If we have a numerical scheme at hand that is capable to treat the source terms properly 
% and we include the sources into the \CK-procedure all the above methods can be applied.
% 
% Especially the source does not affect the decomposition of the temporal derivatives at the junction, 
% since lower order terms in \eqref{eq:generalizedRP} are dropped anyhow.

Source terms can be treated in a straight forward way,
since they do not affect the coupling procedure.
Only when particular steady states should be preserved,
the averaging process in \eqref{eq:LumpingPDE} has to be adapted accordingly.
In case of a balance law we obtain
\begin{align*}
\partial_tu(t,x)+\partial_x\left(f(u(t,x))\right)&=S(u(t,x)) \\
\Rightarrow 
\partial_t
U
+\frac{1}{L}
\left[
f(U_g^r)-f(U_g^l)
\right]&=\frac{1}{L} \int_{0}^{L} S(u(t,x)) \approx \tilde{S}(U) \
\end{align*}
For example in the case of the shallow water equations \eqref{eq:SWE} a huge variety of well-balanced numerical schemes has been developed 
in order to incorporate the bottom elevation in a suitable manner e.g. \cite{Bermudez19941049,audusse_fast_2004,Canestrelli2009834}.
In this particular situation the approximation made to obtain $\tilde{S}$ reads as
\begin{align}\label{eq:approx_bootom}
\frac{1}{L}\int_{0}^{L}-gh\partial_xb&\approx -g \bar{h} \frac{b(L)-b(0)}{L}\ ,
\end{align}
where $h$ denotes the water level, $b$ the bottom elevation and $g$ the gravitational acceleration.
\begin{figure}
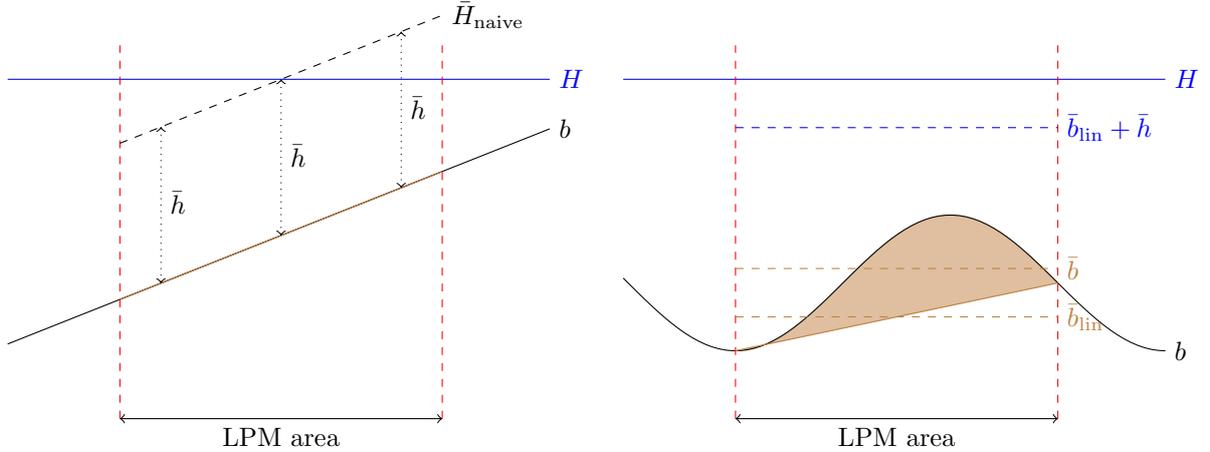
 
	\centering
	\externaltikz{LPM-naive}{
		\relinput{Tikzbilder/Naive-lumping/Naive-lumping.tex}
	}
	\caption{Problems arising from naive lumping. Left: Scheme can become very unbalanced. Right: Conservation of mass violated.}
	\label{fig:naive LPM problems}
\end{figure}

%\begin{figure} 
%	\centering
%	\externaltikz{LPM-not-well-balanced-schematic}{
%		\relinput{Tikzbilder/LPM_well_balancing/LPM_well_balancing.tex}
%	}
%	\caption{Naive treatment leads to a very unbalanced scheme}
%	\label{fig:LPM unbalanced schematic}
%\end{figure}
As illustrated in figure \ref{fig:naive LPM problems} an independent treatment of flux an source
can not lead to a method preserving steady states.
Therefore we have to apply a reconstruction technique to the values $U$, which follows the bottom topography.
This so called hydrostatic reconstruction \cite{audusse_fast_2004} is a well known tool in this context.
Thus whenever bottom elevation is considered we will apply the hydrostatic reconstruction on lumped edges
in order to capture the 'lake at rest' correctly.
% Furthermore, ADER schemes are asymptotically well balanced \todo{Referenz}, 
% the problem therefore disappears with higher resolution, 
% which of course does not apply for the LPM model. 

% \subsubsection{Preventing mass loss due to lumping}
%\begin{figure} 
%	\centering
%	\externaltikz{LPM-mass-loss}{
%		\relinput{Tikzbilder/LPM_bottom_discrepancy/LPM_bottom_discrepancy.tex}
%	}
%	\caption{Mass loss due to naive lumping}
%	\label{fig:LPM mass loss schematic}
%\end{figure}
A further detail concerning this particular steady state is also depicted in figure \ref{fig:naive LPM problems}.
If an initial condition for the PDE model is provided on the network,
the initial conditions of the lumped parameter model are obtained by the averaging process in \eqref{eq:LumpingPDE}.
Since in \eqref{eq:approx_bootom} the bottom elevation is linearized, 
possible errors have to be compensated in the initial values of $U$.
In the following we will just add the missing amount of water artificially to the initial states of the LPM model 
by modifying the reconstruction to
\begin{align*}
	U&=\frac{1}{L}
	\int_{0}^{L} u(x) dx+
	\left(
			\frac{1}{L}\int_{0}^{L} b(x) dx
			-\frac{1}{2}(b(L)+b(0))
	\right) \	.		
\end{align*}

\section{Numerical Examples}
\label{S:Numerik}
In this section we investigate the above presented numerical methods in different test cases.
As conservation law along the edges we choose the shallow water equations 
\begin{align}
 \label{eq:SWE}
% 	\left.
	\begin{aligned}		
	\partial_t h + \partial_x q&=0 \\
	\partial_t q + \partial_x \left(\frac{q^2}{h}+\frac{1}{2}gh^2\right)&=0 \ .
	\end{aligned} 	
% 	\right\rbrace & \text{PDE}	
% 	\\
\end{align}
$h$ denotes the depth of the water and $q$ is the discharge in $x$ direction.
$g=9.81$ is the gravitational acceleration.

As coupling conditions in the junctions we consider two different sets of equations.
The first one is the so called 'equal heights' coupling, which reads for $n$ connected edges
% \begin{Definition}
% 	Pressure coupling %\cite{MR2223073,MR2304344} :%\\
% 	Equal hight coupling %\cite{MR2223073,MR2304344} :%\\
	\begin{align}
% 	\label{def:splitter vertex}
      \label{eq:equalheights}
	\begin{split}
	\sum_{i=1}^nq_i         & =0 \ ,                          \\
% 	p_1(\rho_1)-p_i(\rho_i) & =0  \quad 2\leq i \leq n \ . 
	h_1-h_i & =0  \quad 2\leq i \leq n \ . 
	\end{split}
	\end{align}
% 	In case of an identical pressure law in all connected edges, the last $n-1$ equations reduce to $\rho_1-\rho_i=0  \quad 2\leq i \leq n $.
% \end{Definition}
The first equation states the conservation of mass at the junction. 
The remaining $n-1$ equations force all heights at the junction to be at the same level.

The second set of coupling conditions involve an ODE at the junction.
Here we consider a storage tank or manhole at the coupling point, which is modeled by
\begin{align}
%  	\begin{aligned}
\nonumber
		\dot h_m&=\frac{Q_m}{A_m} \\
		\label{eq:Manhole_Qm}
		\dot Q_m&=\frac{gA_m}{h_m}
		\left(
		\frac{1}{2g}\frac{q_1^2}{h_1^2}+h_1
		-\left(
		\frac{1}{2g}\frac{Q_m^2}{A_m^2}+h_m
		\right)
		\right) \ .
% 	\end{aligned}
\end{align}
The two states are the vertical water level in the tank $h_m$ and the discharge $Q_m$ flowing into the volume.
The constant $A_m$ is the horizontal cross sectional area of the storage tank.
It is important that this model is always accompanied by the following set of coupling conditions
% \begin{Definition}{Energy conserving Manhole model}
% \label{def:full manhole model}
% The manhole model we use as example is defined as follows:\\
% 	\todo{Zitation}
\begin{align*}
% 	\left.
	\begin{aligned}
% 	\sum_{i=1}^{n} q_ib_i+Q_m&=0 \\
	\sum_{i=1}^{n} q_i+Q_m&=0 \ ,\\
	\frac{1}{2g}\frac{q_1^2}{h_1^2}+h_1-\left(\frac{1}{2g}\frac{q_i^2}{h_i^2}+h_i\right)&=0 
	  \qquad 2\leq i \leq n \ .
% 	\\
	\end{aligned}
% 	\right\rbrace & \text{ACC}	 
% 	\\
% 	\left.
% 	\begin{aligned}
% 		\dot h_m&=\frac{Q_m}{A_m} \\
% 		\dot Q_m&=\frac{gA_m}{h_m}
% 		\left(
% 		\frac{1}{2g}\frac{q_i^2}{h_i^2}+h_i
% 		-\left(
% 		\frac{1}{2g}\frac{Q_m^2}{A_m^2}+h_m
% 		\right)
% 		\right) 
% 	\end{aligned}
% 	\right\rbrace & \text{ODE}	 
\end{align*}
% \end{Definition}
The first equation again ensures the conservation of total mass in the coupled system.
The following $n-1$ equations state the equality of the so called hydraulic heads or energy levels.
As this name indicates these conditions provide the conservation of the total energy at the junction in case of smooth solutions \cite{MR3335526,MR3281585}.
The total energy in the coupled system is conserved due to equation \eqref{eq:Manhole_Qm} in the storage model \cite{MR3281585}.
Here we further note that \eqref{eq:Manhole_Qm} does not depend on the choice of the related edge 
since all hydraulic heads coincide.

\begin{figure}[]
	\centering
	\externaltikz{split-circle-network}{
		\input{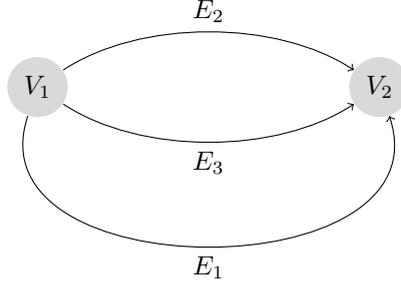}
	}
	\caption{Split circle network.}
	\label{fig:splitted circle network}
\end{figure}
In the following convergence studies, we use a simple network consisting of three edges and two nodes,
as shown in \figref{fig:splitted circle network}.
All edges are of the same length $L=25$ 
and in both nodes a storage tank with $A_m=1$ is placed.
The initial data used for convergence studies needs not only to be smooth along the edges, 
but also has to satisfy the coupling conditions and its temporal derivatives up to the order of the schemes to be investigated.
In the following we choose lake at rest like states at the junctions
and a sufficiently smooth transition between their two water levels.
This can be achieved by $q_i(x,0)=0$ $i=1,2,3$ 
and polynomial states  of degree $15$ for $h_i$ which are determined by the constraints
% with suitably smooth initial data on edges of length $L=25$. 
% The only initial data we could find
% that fulfils the coupling conditions up to arbitrary order is water at rest, 
% therefore we use initial data defined by polynomials of sufficient degree with:
\def\linkerwert{2}
\def\rechterwert{3}
\def\tende{2.4}
\begin{align*}
% h_i(0,0)&=\linkerwert &h_i(L,0)&=\rechterwert & q_i(x,0)&=0\\
% \partial_x^kh_i(0,0)&=0 & \partial_x^kh_i(L,0)&=0 &
h_i(0,0)=\linkerwert\ , \quad \partial_x^kh_i(0,0)=0\ ,  
&&
h_i(L,0)=\rechterwert\ , \quad \partial_x^kh_i(L,0)=0
\ , \quad k = 1,\dots,7\ . 
\end{align*}
The initial states of the ODEs are also at rest, i.e. $w^1(0)=(2,0)^T$ and $w^2(0)=(3,0)^T$.
\begin{figure}
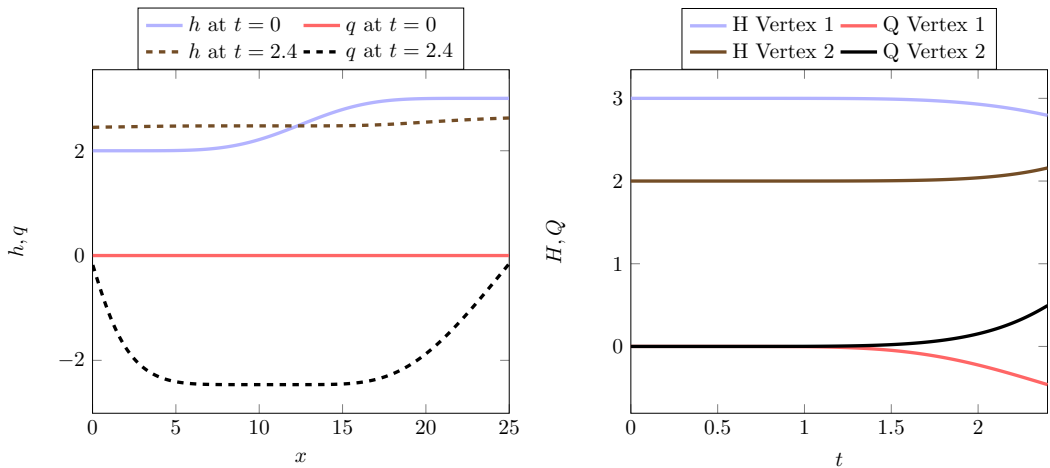

	\centering
	\externaltikz{Data-and-solution-convergence-study-splitcircle}{
	\begin{tikzpicture}[scale=0.8]% file
				\def\axisarguments{name=axis1};
			    \input{Tikzbilder/initial-data-convergence-study.tex}			   				
				\def\axisarguments{name=axis2,at={(axis1.east)},anchor=west,xshift=2cm,};
				\input{Tikzbilder/ODE-data-convergence-study.tex}					
	\end{tikzpicture}
	}
	\caption{Initial Data and solution for the convergence studies \ref{sec:Convergence study TT} and \ref{sec:Convergence study HEOC}. PDE states on the edges to the left, ODE states in the vertices on the right side.} 
	\label{fig:Data and solution split circle convergence study}
\end{figure}
 \MyFloatBarrier	
In \figref{fig:Data and solution split circle convergence study} we show initial states
as well as the reference solution for the PDEs on the left side,
on the right hand side are the states of the ODEs in the nodes for $t\in[0,\tende]$.
The computation of the reference solution was performed by a scheme of order $8$ on a grid of $800$ cells per edge.

In all numerical examples the time step $\Delta t$ is synchronized in the complete network according to a CFL number $0.95$.
The stability bounds of the ODEs are always less restrictive than those of the PDEs.

\subsection{Convergence study \TT}
\label{sec:Convergence study TT}
\begin{figure}
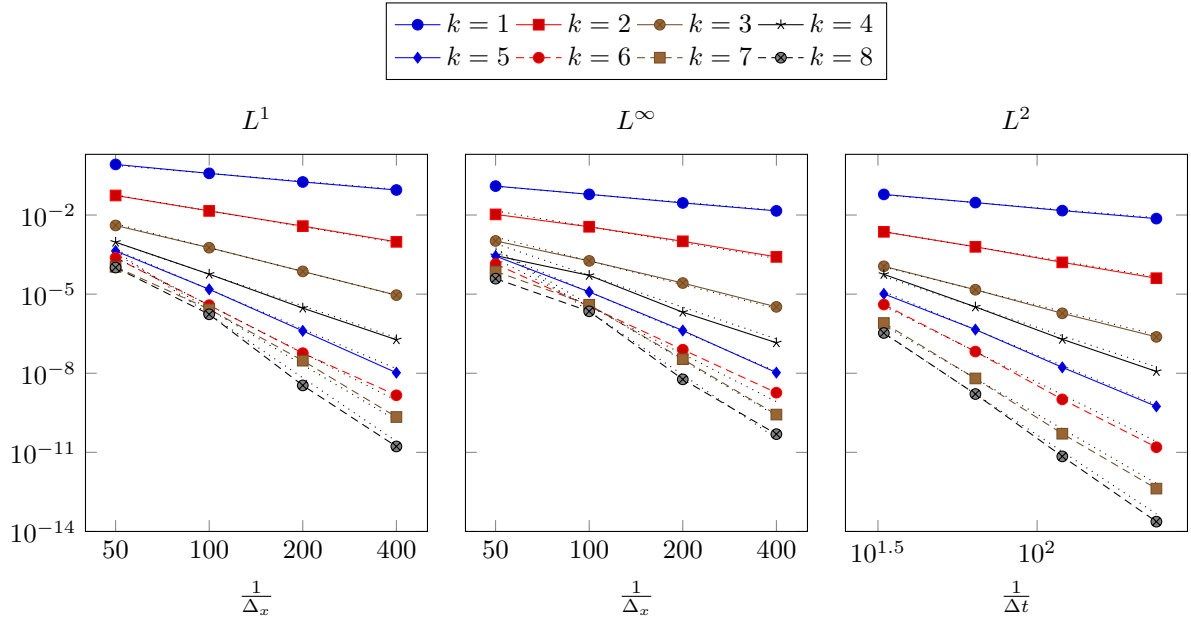

 	\centering
 	\externaltikz{Konvergenzplot-splitcircle-TT-klein}{
 		\begin{tikzpicture}[scale=1]
 		\def\zusatz{\legend{}};
 		\def\gemeinsameargumente{width=4.5cm,height=5cm,ymin=0.00000000000001,ymax=2,ylabel={}};
 		\def\axisarguments{name=axis1,title=$L^1$,\gemeinsameargumente};
 		\input{Tikzbilder/Konvergenztest-1norm.tex};			   				
 		\def\axisarguments{name=axis2,at={(axis1.east)},anchor=west,xshift=0.5cm,yticklabels={,,},\gemeinsameargumente};
 		\def\zusatz{};
 		\input{Tikzbilder/Konvergenztest-infnorm.tex};
 		\def\axisarguments{name=axis3,at={(axis2.east)},anchor=west,xshift=0.5cm,yticklabels={,,},title=$L^2$,\gemeinsameargumente};
 		\def\zusatz{\legend{}};
 		\input{Tikzbilder/Konvergenztest-ODE}		
 		\end{tikzpicture}
 	}
 	\caption{Convergence study \ref{sec:Convergence study TT}: Splitcircle TT solver in $L^1$ and $L^\infty$ norm for the edges, $L^2$ norm for the ODEs.}
 	\label{fig:convergence plots splitcircle TT}
 \end{figure}

 The first convergence tests we perform for the \TT solver which is described in section \ref{S:GRP-TT}.
 In \figref{fig:convergence plots splitcircle TT} the $L^1$ and the $L^\infty$-error at $t=\tende$ are plotted against the reciprocal of the cell width $\Delta x$.
 Furthermore the $L^2$-errors of the involved ODEs are shown. 
  %\tikzremakenext
 Additionally for selected orders $k$ the errors and rates of convergence are given in \tabref{tab:convergence table splitcircle TT}.  
 \begin{table}
 \centering
 \begin{tabular}{|c|c|c||c|c||c|c||c|c||c|c||}
 	\cline{2-11}
 	%\hline
 	\multicolumn{1}{c}{} & \multicolumn{2}{|c||}{$k=2$} & \multicolumn{2}{|c||}{$k=4$} & \multicolumn{2}{|c||}{$k=5$} & \multicolumn{2}{|c||}{$k=6$} & \multicolumn{2}{|c||}{$k=8$} \\
 	\cline{2-11}
 	\multicolumn{8}{l}{PDE $L^1$ norm} \\
 	\hline
 	N                    & $L^1$                        & $O_{L^1}$                    & $L^1$                        & $O_{L^1}$                    & $L^1$                        & $O_{L^1}$      & $L^1$      & $O_{L^1}$      & $L^1$      & $O_{L^1}$      \\
 	\hline
 	50                   & 5.51e-02                     &                              & 9.32e-04                     &                              & 4.34e-04                     &                & 2.37e-04   &                & 1.02e-04   &                \\
 	\hline
 	100                  & 1.44e-02                     & 1.94                         & 5.71e-05                     & 4.03                         & 1.52e-05                     & 4.83           & 3.81e-06   & 5.96           & 1.71e-06   & 5.89           \\
 	\hline
 	200                  & 3.79e-03                     & 1.92                         & 2.95e-06                     & 4.28                         & 4.04e-07                     & 5.24           & 5.77e-08   & 6.04           & 3.45e-09   & 8.96           \\
 	\hline
 	400                  & 9.61e-04                     & 1.98                         & 1.86e-07                     & 3.98                         & 1.05e-08                     & 5.27           & 1.46e-09   & 5.31           & 1.67e-11   & 7.69           \\
 	\hline
 %	\multicolumn{1}{c}{} & \multicolumn{2}{|c||}{$k=2$} & \multicolumn{2}{|c||}{$k=4$} & \multicolumn{2}{|c||}{$k=5$} & \multicolumn{2}{|c||}{$k=6$} & \multicolumn{2}{|c||}{$k=8$}  \\
 	\multicolumn{8}{l}{PDE $L^\infty$ norm} \\
 	\hline
 	N                    & $L^\infty$                   & $O_{L^\infty}$               & $L^\infty$                   & $O_{L^\infty}$               & $L^\infty$                   & $O_{L^\infty}$ & $L^\infty$ & $O_{L^\infty}$ & $L^\infty$ & $O_{L^\infty}$ \\
 	\hline
 	50                   & 1.06e-02                     &                              & 2.89e-04                     &                              & 2.77e-04                     &                & 1.38e-04   &                & 3.85e-05   &                \\
 	\hline
 	100                  & 3.59e-03                     & 1.56                         & 5.13e-05                     & 2.50                         & 1.21e-05                     & 4.51           & 3.43e-06   & 5.34           & 2.23e-06   & 4.11           \\
 	\hline
 	200                  & 1.01e-03                     & 1.83                         & 2.06e-06                     & 4.64                         & 4.17e-07                     & 4.86           & 7.81e-08   & 5.45           & 5.85e-09   & 8.57           \\
 	\hline
 	400                  & 2.58e-04                     & 1.96                         & 1.43e-07                     & 3.85                         & 1.05e-08                     & 5.30           & 1.83e-09   & 5.42           & 4.90e-11   & 6.90           \\
 	\hline
 	\multicolumn{8}{l}{ODE $L^2$ norm} \\
 	\hline
 	N                    & $L^2$                        & $O_{L^2}$                    & $L^2$                        & $O_{L^2}$                    & $L^2$                        & $O_{L^2}$      & $L^2$      & $O_{L^2}$      & $L^2$      & $O_{L^2}$      \\
 	\hline
 	33                   & 2.31e-03                     &                              & 5.62e-05                     &                              & 1.02e-05                     &                & 4.02e-06   &                & 3.35e-07   &                \\
 	\hline
 	64                   & 6.22e-04                     & 1.98                         & 3.27e-06                     & 4.30                         & 4.53e-07                     & 4.70           & 6.64e-08   & 6.19           & 1.63e-09   & 8.05           \\
 	\hline
 	120                  & 1.60e-04                     & 2.16                         & 1.94e-07                     & 4.50                         & 1.67e-08                     & 5.25           & 1.02e-09   & 6.64           & 7.02e-12   & 8.66           \\
 	\hline
 	237                  & 4.01e-05                     & 2.04                         & 1.17e-08                     & 4.13                         & 5.53e-10                     & 5.01           & 1.56e-11   & 6.14           & 2.35e-14   & 8.37           \\
 	\hline
 \end{tabular}
 \caption{Convergence study \ref{sec:Convergence study TT}: Convergence rates splitcircle with TT solver.}
 \label{tab:convergence table splitcircle TT}
 \end{table}
% To demonstrate the convergence of the PDEs, we plot the error against the reciproke of the cell width 
%  with auxiliary lines indicating the expected slope for $L^1$ and $L^\infty$ norms (see \figref{fig:convergence plots splitcircle TT}).
 All numerical solutions converge with the expected order.
%  To check the convergence of the ODEs in the vertices, 
%  we plot the error at $t=\tende$ against $\frac{1}{\Delta t}$.
%  Again we can obsorve the expected order of convergence, 
%  see \figref{fig:convergence plots splitcircle TT} 
%  and \tabref{tab:convergence table splitcircle TT}.

% \tikzremakenext
% \begin{figure}
% 	\centering
% 	\externaltikz{Konvergenztest-ODE}{
% 		\begin{tikzpicture}
% 		\def\zusatz{}
% 		\def\axisarguments{}
% 		\input{Tikzbilder/Konvergenztest-ODE}
% 		\end{tikzpicture}
% 	}
% 	\caption{convergence rates TT ODE}
% 	\label{fig:convergence plot TT ODE}
% \end{figure}
 
 \MyFloatBarrier
 \subsection{Convergence study \HEOC}
 \label{sec:Convergence study HEOC}
 For the HEOC-coupling procedure we repeat the same test for schemes up to order $k=6$.
 The $L^1$- and $L^\infty$ errors for the PDEs and the $L^2$-errors of the ODEs are shown in in \figref{fig:convergence plots splitcircle HEOC}.
 The precise values and the corresponding convergence rates are presented in \tabref{tab:convergence table splitcircle HEOC}.
 \begin{figure}
 	\centering
 	\externaltikz{Konvergenzplots-splitcircle-HEOC-klein}{
 		\begin{tikzpicture}
 		\def\zusatz{\legend{}};
 		\def\gemeinsameargumente{width=4.5cm,height=5cm,ymin=0.000000000001,ymax=2,ylabel={}};
 		\def\axisarguments{name=axis1,title=$L^1$,\gemeinsameargumente};
 		\input{Tikzbilder/Konvergenztest-HEOC-1norm}			   				
 		\def\axisarguments{name=axis2,at={(axis1.east)},anchor=west,xshift=0.5cm,yticklabels={,,},\gemeinsameargumente};
 		\def\zusatz{};
 		\input{Tikzbilder/Konvergenztest-HEOC-infnorm.tex}
 		\def\axisarguments{name=axis3,at={(axis2.east)},anchor=west,xshift=0.5cm,yticklabels={,,},title=$L^2$,\gemeinsameargumente};
 		\def\zusatz{\legend{}};
 		\input{Tikzbilder/Konvergenztest-HEOC-ODE.tex}		
 		\end{tikzpicture}		
 	}
\caption{Convergence study \ref{sec:Convergence study HEOC}: Splitcircle HEOC solver in $L^1$ and $L^\infty$ norm for the edges, $L^2$ norm for the ODEs.} 	
 	\label{fig:convergence plots splitcircle HEOC}
 \end{figure}
 \begin{table}
 	\centering
 \begin{tabular}{|c||c|c||c|c||c|c||c|c||} 
 	\cline{2-9}
 	\multicolumn{1}{c}{}  &\multicolumn{2}{|c||}{$k=2$}  &\multicolumn{2}{|c||}{$k=4$}  &\multicolumn{2}{|c||}{$k=5$}  &\multicolumn{2}{|c||}{$k=6$}  \\ 
 	\cline{2-9}
 	\multicolumn{9}{l}{PDE $L^1$ norm} \\ 
 	\hline 
 	N&$L^1$ & $O_{L^1}$&$L^1$ & $O_{L^1}$&$L^1$ & $O_{L^1}$&$L^1$ & $O_{L^1}$\\ 
 	
 	\hline
 	50                   & 4.39e-02                     &                              & 6.80e-04                     &                              & 1.38e-04 &           & 1.94e-04 &           \\
 	\hline
 	100                  & 1.13e-02                     & 1.95                         & 5.21e-05                     & 3.71                         & 5.09e-06 & 4.76      & 3.48e-06 & 5.80      \\
 	\hline
 	200                  & 3.03e-03                     & 1.90                         & 2.97e-06                     & 4.13                         & 1.48e-07 & 5.10      & 1.74e-08 & 7.64      \\
 	\hline
 	400                  & 7.69e-04                     & 1.98                         & 2.02e-07                     & 3.88                         & 3.63e-09 & 5.35      & 3.18e-10 & 5.78      \\
 	\hline
 	
 	\hline 
 	\multicolumn{9}{l}{PDE $L^\infty$ norm} \\ 
 	\hline 
 	N&$L^\infty$ & $O_{L^\infty}$&$L^\infty$ & $O_{L^\infty}$&$L^\infty$ & $O_{L^\infty}$&$L^\infty$ & $O_{L^\infty}$\\ 
 	
 	\hline
 	50                   & 1.05e-02                     &                              & 3.49e-04                     &                              & 6.87e-05   &                & 1.40e-04   &                \\
 	\hline
 	100                  & 3.22e-03                     & 1.70                         & 4.59e-05                     & 2.93                         & 5.97e-06   & 3.52           & 9.02e-06   & 3.96           \\
 	\hline
 	200                  & 8.73e-04                     & 1.88                         & 1.62e-06                     & 4.83                         & 1.77e-07   & 5.08           & 2.65e-08   & 8.41           \\
 	\hline
 	400                  & 2.22e-04                     & 1.98                         & 1.20e-07                     & 3.76                         & 3.52e-09   & 5.65           & 1.50e-09   & 4.14           \\
 	\hline
 	
 	\hline 
 	\multicolumn{9}{l}{ODE $L^2$ norm} \\ 
 	\hline 
 	N&$L^2$ & $O_{L^2}$&$L^2$ & $O_{L^2}$&$L^2$ & $O_{L^2}$&$L^2$ & $O_{L^2}$\\ 
 	
 	\hline
 	33                   & 2.86e-03                     &                              & 1.06e-05                     &                              & 7.86e-06 &           & 6.32e-06 &           \\
 	\hline
 	64                   & 8.16e-04                     & 1.89                         & 1.25e-06                     & 3.24                         & 2.79e-07 & 5.04      & 3.85e-07 & 4.22      \\
 	\hline
 	120                  & 2.13e-04                     & 2.13                         & 2.90e-07                     & 2.32                         & 9.12e-09 & 5.44      & 1.23e-09 & 9.14      \\
 	\hline
 	237                  & 5.40e-05                     & 2.02                         & 2.08e-08                     & 3.87                         & 2.63e-10 & 5.21      & 8.86e-12 & 7.25      \\
 	\hline
 	
 	\hline 
 \end{tabular}
 \caption{Convergence study \ref{sec:Convergence study HEOC}: Convergence rates Splitcircle HEOC solver.}
 \label{tab:convergence table splitcircle HEOC}
 \end{table} 
% \begin{figure}
% 	\centering
% 	\externaltikz{Konvergenzplots-splitcircle-HEOC-gross}{
% 	\begin{tikzpicture}[scale=0.8]
% 	\def\zusatz{}
% 	\def\axisarguments{name=axis1};
% 	\input{Tikzbilder/Konvergenztest-HEOC-1norm}			   				
% 	\def\axisarguments{name=axis2,at={(axis1.east)},anchor=west,xshift=2cm,};
% 	\input{Tikzbilder/Konvergenztest-HEOC-infnorm.tex}
% 	\def\axisarguments{name=axis3,at={(axis1.south)},anchor=north,yshift=-3cm};
% 	\input{Tikzbilder/Konvergenztest-HEOC-ODE.tex}		
% 	\end{tikzpicture}
% 	}
% 	\caption{Convergence plots splitcircle HEOC in $L^1$ and $L^\infty$ norm for the edges, $L^2$ norm for the ODEs }
% 	\label{fig:convergence plots splitcircle HEOC}
% \end{figure}

The errors are within the same range as those of the \TT method
and all solutions converge with the predicted order.
 
 \MyFloatBarrier

 \subsection{Convergence study Lumped Parameter Model}
 \label{sec:Convergence study LPM}
 The final convergence test addresses lumped parameter models.
 \begin{figure}
 	\centering
 \externaltikz{Diamond-LPM-skizze}{
 \begin{tikzpicture}[scale=2]
 \tikzstyle{knoten}=[circle, fill=gray!30, scale=0.7];
 \tikzstyle{kante}=[-,line width=1];
 \node[knoten,fill=blue,opacity=0.4] (1) at (2.31,-0.4) {$V_{1}$};
 \node[knoten] (2) at (-0.68,-0.34) {$V_{2}$};
 \node[knoten,fill=blue,opacity=0.4] (3) at (0.72,0.82) {$V_{3}$};
 \node[knoten] (4) at (0.75,-1.09) {$V_{4}$};
 \draw[kante] (1) .. controls(1.8245,-2.5531)and(-0.38,-2.22).. (2) node[pos=0.5,below]{$E_1$};
 \draw[kante] (3) .. controls(1.9,-0.01)and(1.89,0.45).. (1) node[pos=0.5,above right]{$E_2$};
 \draw[kante] (1) .. controls(1.25,-0.78)and(1.08,-1.03).. (4) node[pos=0.5,below right]{$E_3$};
 \draw[kante] (4) .. controls(0.33,0.03)and(0.7,-0.37).. (3) node[pos=0.5,left]{$E_4$};
 \draw[kante] (2) .. controls(-0.11,0.15)and(-0.76,0.05).. (3) node[pos=0.5,left]{$E_5$};
 \draw[kante] (2) .. controls(-0.16,-0.83)and(-0.61,-0.6).. (4) node[pos=0.5,above right]{$E_6$};
 \draw[fill=green,opacity=0.2] (0.8cm,0) ellipse (2cm and 1.3cm);
 \end{tikzpicture}
 }
 \caption{Diamond network. The blue vertices $V_1$ and $V_3$ are energy conserving manholes, i.e. ODE vertices, the rest algebraic equal height coupling vertices. The green area is lumped resulting in a network of one remaining Edge $E_1$ and an LPM vertex with a state of 14 components.}
 \label{fig:Diamond LPM network}
 \end{figure}
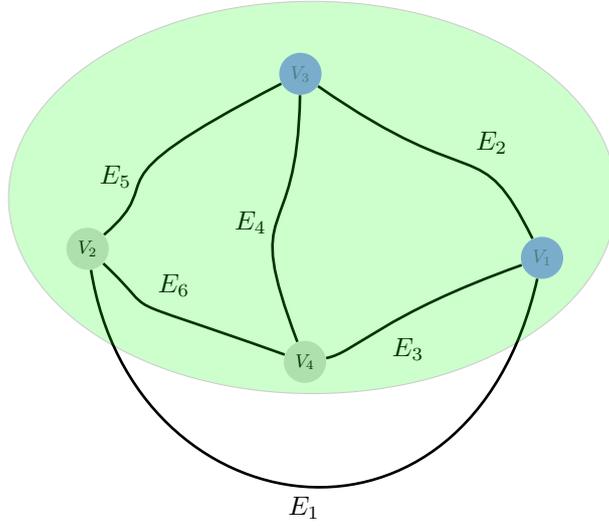
 Therefore we consider a network consisting of six edges and four vertices as depicted in \figref{fig:Diamond LPM network}.
 In the two junctions $V_2$ and $V_4$ equal height coupling is applied,
 while in $V_1$ and $V_3$ two manholes with $A_m=1$ are located.
 The lumping of section \ref{S:Lumped Parameter Models} is applied to the whole network except $E_1$,
 i.e. the green colored region.
 
 \begin{figure}
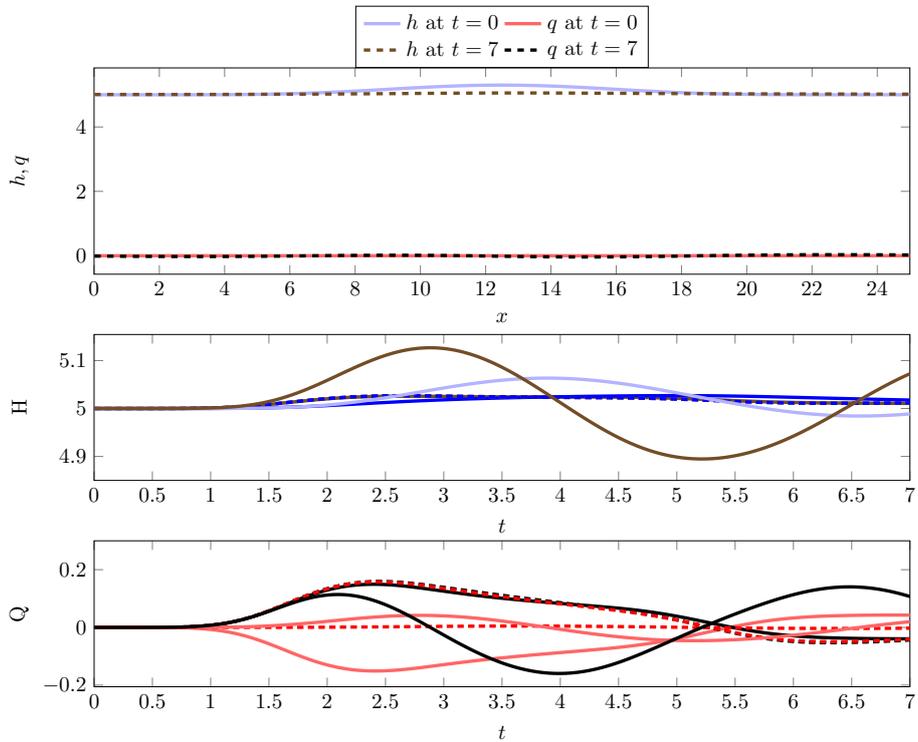

 	\centering
 	\externaltikz{Diamond-LPM-initial-data}{
 		\begin{tikzpicture}[scale=0.8]% file
 		\def\axisarguments{name=axis1,height=5cm,width=15cm};
 		\input{Tikzbilder/Konvergenztest-Diamond-LPM/initial-data-final-data.tex}		   				
 		\def\axisarguments{name=axis2,at={(axis1.south)},anchor=north,yshift=-1cm,ymin=4.85,ylabel=H,height=4cm,width=15cm,};
 		\input{Tikzbilder/Konvergenztest-Diamond-LPM/ODE-data-Diamond-LPM.tex};	
 		\def\axisarguments{name=axis3,at={(axis2.south)},anchor=north,yshift=-1cm,ymax=0.3,ylabel=Q,height=4cm,width=15cm,};
 		\input{Tikzbilder/Konvergenztest-Diamond-LPM/ODE-data-Diamond-LPM.tex};			
 		\end{tikzpicture}
 	}
 	\caption{Convergence study Diamond LPM \ref{sec:Convergence study LPM}, initial data and solutions: From top to bottom: PDE initial data on $E_1$, ODE components of height type and ODE components of impulse type.} 
 	\label{fig:Data and solution Diamond convergence study}
 \end{figure} 
 
 As initial conditions we choose water at rest $q_i(x,0)\equiv 0$ $i=1,\dots,6$
 with the constant water level $h_i(x,0)\equiv 0$, $i = 2,\dots,6$ ,
except for the first edge. 
There we take as data a polynomial of degree $14$ such that
 \begin{align*}
 h_1(0,0)&=5\ ,  && h_1(\frac{L}{2},0)=5.3\ , &&& h_1(L,0)&=5\ , \\
 \partial_x^kh_1(0,0)&=0 \ ,&&&&&\partial_x^kh_1(L,0)&=0 \quad \forall k=1,\dots,6
 \end{align*}
 hold.
 In \figref{fig:Data and solution Diamond convergence study} initial data 
 and solutions are shown.
 For edge $E_1$ only the solution at $t=0$ and $t=7$ are plotted,
 whereas the states of all ODEs are shown on the full time interval $[0,7]$.

 %\tikzremakenext
 \begin{figure}
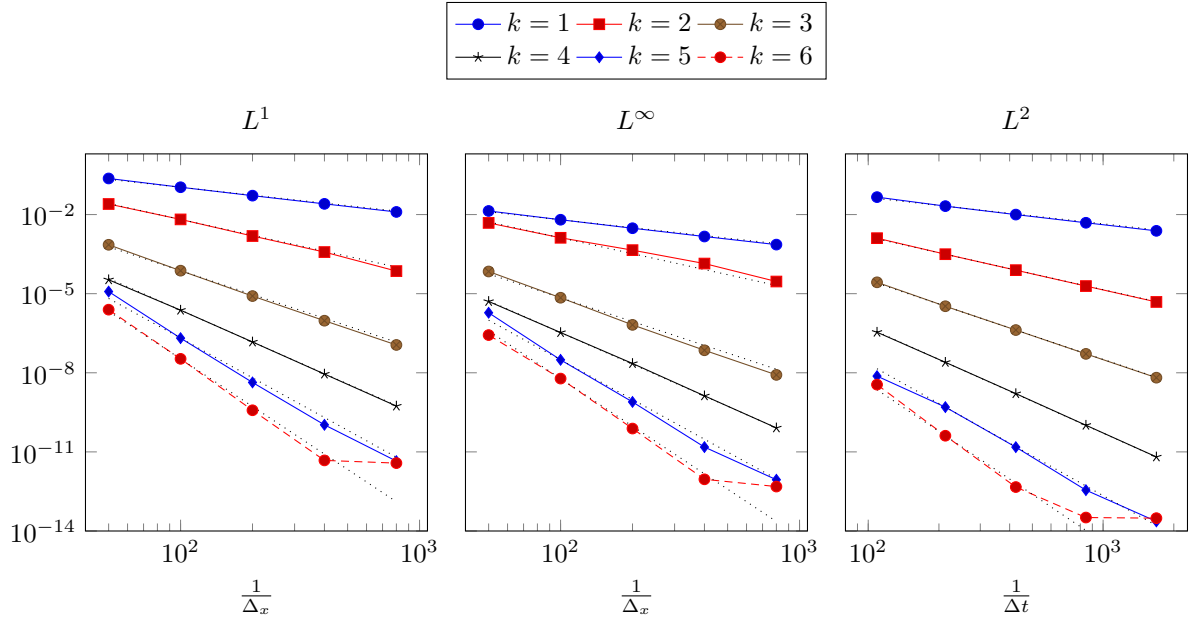

 	\centering
 	\externaltikz{Konvergenzplots-Diamond-LPM-klein}{
 		\begin{tikzpicture}
 		\def\zusatz{\legend{}};
 		\def\gemeinsameargumente{width=4.5cm,height=5cm,ymin=0.00000000000001,ymax=2,ylabel={}};
 		\def\axisarguments{name=axis1,title=$L^1$,\gemeinsameargumente};
 		\input{Tikzbilder/Konvergenztest-Diamond-LPM/Konvergenztest-1norm}			   				
 		\def\axisarguments{name=axis2,at={(axis1.east)},anchor=west,xshift=0.5cm,yticklabels={,,},\gemeinsameargumente};
 		\def\zusatz{};
 		\input{Tikzbilder/Konvergenztest-Diamond-LPM/Konvergenztest-infnorm};
 		\def\axisarguments{name=axis3,at={(axis2.east)},anchor=west,xshift=0.5cm,yticklabels={,,},title=$L^2$,\gemeinsameargumente};
 		\def\zusatz{\legend{}};
 		\input{Tikzbilder/Konvergenztest-Diamond-LPM/Konvergenztest-ODE};
 		\end{tikzpicture}
 	}
 	\caption{Convergence study diamond LPM \ref{sec:Convergence study LPM}, $L^1$ and $L^\infty$ norm for the edge, $L^2$ norm for the ODE.}
 	\label{fig:convergence plots Diamond LPM}
 \end{figure} 
 \begin{table}
 \centering
 \begin{tabular}{|c||c|c||c|c||c|c||c|c||} 
 	\cline{2-9}
 	\multicolumn{1}{c}{}  &\multicolumn{2}{|c||}{$k=2$}  &\multicolumn{2}{|c||}{$k=4$}  &\multicolumn{2}{|c||}{$k=5$}  &\multicolumn{2}{|c||}{$k=6$}  \\ 
 	\cline{2-9}
 	\multicolumn{9}{l}{PDE $L^1$ norm} \\ 
 	\hline 
 	N&$L^1$ & $O_{L^1}$&$L^1$ & $O_{L^1}$&$L^1$ & $O_{L^1}$&$L^1$ & $O_{L^1}$\\ 
 	
 	\hline
 	50                   & 2.51e-02                     &                              & 3.37e-05                     &                              & 1.20e-05 &           & 2.47e-06 &           \\
 	\hline
 	100                  & 6.56e-03                     & 1.94                         & 2.37e-06                     & 3.83                         & 2.06e-07 & 5.87      & 3.38e-08 & 6.19      \\
 	\hline
 	200                  & 1.54e-03                     & 2.09                         & 1.47e-07                     & 4.00                         & 4.31e-09 & 5.58      & 3.79e-10 & 6.48      \\
 	\hline
 	400                  & 3.79e-04                     & 2.02                         & 8.95e-09                     & 4.04                         & 1.06e-10 & 5.34      & 4.76e-12 & 6.32      \\
 	\hline
 	800                  & 7.19e-05                     & 2.40                         & 5.48e-10                     & 4.03                         & 4.47e-12 & 4.57      & 3.75e-12 & 0.34      \\
 	\hline
 	
 	\hline 
 	\multicolumn{9}{l}{PDE $L^\infty$ norm} \\ 
 	\hline 
 	N&$L^\infty$ & $O_{L^\infty}$&$L^\infty$ & $O_{L^\infty}$&$L^\infty$ & $O_{L^\infty}$&$L^\infty$ & $O_{L^\infty}$\\ 
 	
 	\hline
 	50                   & 4.82e-03                     &                              & 5.07e-06                     &                              & 1.89e-06   &                & 2.69e-07   &                \\
 	\hline
 	100                  & 1.32e-03                     & 1.87                         & 3.36e-07                     & 3.91                         & 3.09e-08   & 5.93           & 6.03e-09   & 5.48           \\
 	\hline
 	200                  & 4.47e-04                     & 1.56                         & 2.25e-08                     & 3.90                         & 7.85e-10   & 5.30           & 7.67e-11   & 6.30           \\
 	\hline
 	400                  & 1.39e-04                     & 1.69                         & 1.34e-09                     & 4.08                         & 1.52e-11   & 5.69           & 9.12e-13   & 6.39           \\
 	\hline
 	800                  & 2.89e-05                     & 2.26                         & 8.13e-11                     & 4.04                         & 8.79e-13   & 4.12           & 4.85e-13   & 0.91           \\
 	\hline
 	
 	\hline 
 	\multicolumn{9}{l}{ODE $L^2$ norm} \\ 
 	\hline 
 	N&$L^2$ & $O_{L^2}$&$L^2$ & $O_{L^2}$&$L^2$ & $O_{L^2}$&$L^2$ & $O_{L^2}$\\ 
 	
 	\hline
 	109                  & 1.27e-03                     &                              & 3.45e-07                     &                              & 7.49e-09 &           & 3.56e-09 &           \\
 	\hline
 	213                  & 3.14e-04                     & 2.09                         & 2.50e-08                     & 3.92                         & 5.00e-10 & 4.04      & 4.10e-11 & 6.66      \\
 	\hline
 	425                  & 7.80e-05                     & 2.02                         & 1.62e-09                     & 3.96                         & 1.52e-11 & 5.06      & 4.64e-13 & 6.49      \\
 	\hline
 	845                  & 1.94e-05                     & 2.02                         & 1.03e-10                     & 4.02                         & 3.55e-13 & 5.47      & 3.24e-14 & 3.87      \\
 	\hline
 	1686                 & 4.85e-06                     & 2.01                         & 6.45e-12                     & 4.00                         & 2.31e-14 & 3.95      & 3.07e-14 & 0.08      \\
 	\hline
 	
 	\hline 
 \end{tabular}
 	\caption{Convergence study diamond LPM \ref{sec:Convergence study LPM}: convergence rates.}
 	\label{tab:convergence table Diamond LPM}
 \end{table}
% \begin{figure}
% 	\centering
% 	\externaltikz{Konvergenzplots-Diamond-LPM}{
% 		\begin{tikzpicture}[scale=0.8]
% 		\def\axisarguments{name=axis1};
% 		\input{Tikzbilder/Konvergenztest-Diamond-LPM/Konvergenztest-1norm}			   				
% 		\def\axisarguments{name=axis2,at={(axis1.east)},anchor=west,xshift=2cm,};
% 		\input{Tikzbilder/Konvergenztest-Diamond-LPM/Konvergenztest-infnorm};			 
% 		\def\axisarguments{name=axis3,at={(axis1.south)},anchor=north,yshift=-3cm};
% 		\input{Tikzbilder/Konvergenztest-Diamond-LPM/Konvergenztest-ODE};	 
% 		\end{tikzpicture}
% 	}
% 	\caption{Convergence plots Diamond LPM in $L^1$ and $L^\infty$ norm for the edge, $L^2$ norm for the ODE }
% 	\label{fig:convergence plots Diamond LPM}
% \end{figure} 
  The error plots are shown in \figref{fig:convergence plots Diamond LPM},
  with the corresponding data in \tabref{tab:convergence table Diamond LPM}.
  The solutions in all components involved converge with the designed order.
 \MyFloatBarrier
  
 \subsection{Capturing of shock waves}
 \label{sec:capturing shock waves}
 In regions of smooth states the advantages of higher order methods are clearly indicated by the order of convergence.
This does not hold for discontinuous solutions.
In the following example we want to investigate the stability and accuracy of high order methods near shock waves.

 \begin{figure}
 	\centering
 \externaltikz{Modified-splitcircle}{
 \begin{tikzpicture}
 \tikzstyle{knoten}=[circle, fill=gray!30, scale=1];
 \tikzstyle{kante}=[-,line width=1];
 \node[knoten] (1) at (-2,0) {$V_{1}$};
 \node[knoten] (2) at (2,0) {$V_{2}$};
 \node[knoten] (3) at (-1.01,-2.99) {$V_{3}$};
 \draw[kante] (1) .. controls(-0.66,-0.87)and(0.75,-0.69).. (2) node[pos=0.5,above]{$E_1$};
 \draw[kante] (1) .. controls(-0.74,0.81)and(0.76,0.67).. (2) node[pos=0.5,above]{$E_2$};
 \draw[kante] (3) .. controls(-1.82,-2.07)and(-2.05,-1.18).. (1) node[pos=0.5,left]{$E_3$};
 \draw[kante] (2) .. controls(2.84,-2.64)and(1.81,-3.74).. (3) node[pos=0.5,above]{$E_4$};
 \end{tikzpicture}
 }
 \caption{Modified splitcircle.}
 \label{fig:modified splitcircle}
 \end{figure}
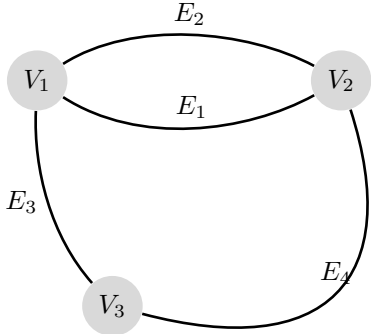 
Therefore we consider a modified split circle network as depicted in \figref{fig:modified splitcircle}.
In all the nodes the coupling conditions for a storage tank with $A_m=1$ are applied. 
%  We want to demonstrate the advantages of high order schemes even in presence of shocks 
%  which limit the order of convergence to one.
%  Therefore we choose a modified splitcircle \figref{fig:modified splitcircle}
%  with edges governed by shallow water equations and manhole vertices \ref{def:full manhole model}.
As initial conditions we choose the constant water levels
 With the following initial data:
 \begin{align*}
 h^1(x,0)&=h^2(x,0)=h^4(x,0)=5 \ ,&
 h^3(x,0)&=6 
%  \\
%  q^1(0,x)&=q^2(0,x)=q^3(0,x)=q^4(0,x)=0
 \end{align*} 
 and $q^i(x,0) \equiv 0$ $ i = 1,\dots,4$.
 Using this setup instead of reusing the normal split circle with Riemann initial data ensures a clean Riemann problem at the junction without the interference of intermediary states. 
  
\begin{figure}
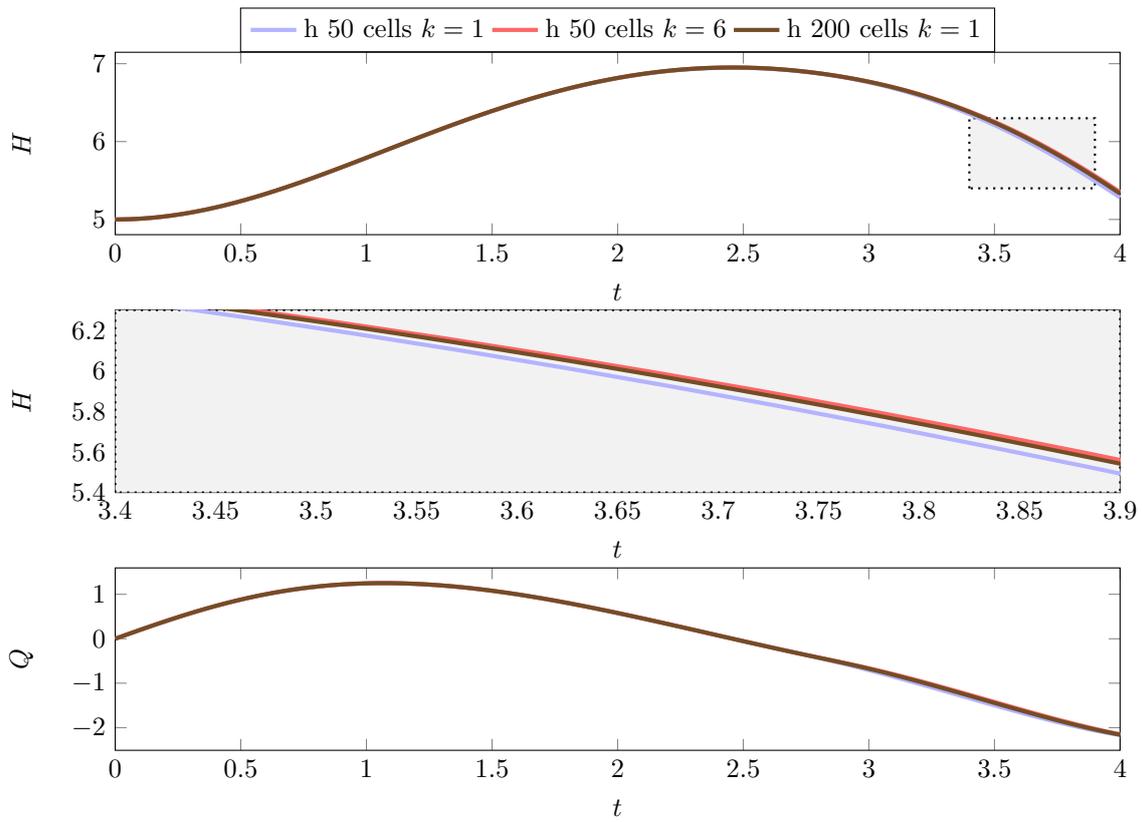
	
 	\centering		
 	\externaltikz{Splitcircle-shock-V1}{
 		\def\ausschnittxmin{3.4}
 		\def\ausschnittxmax{3.9}
 		\def\ausschnittymin{5.4}
 		\def\ausschnittymax{6.3}
 		\begin{tikzpicture}
 		\def\axisarguments{name=A1}
 		\def\precode{
 			\draw[thick,fill=gray!10,style=dotted] (axis cs: \ausschnittxmin,\ausschnittymin) rectangle (axis cs: \ausschnittxmax,\ausschnittymax);
 			}
 		\def\postcode{}
 		\input{Tikzbilder/Splitcircle-Schock/V1-plot-H.tex}
 		% % zweites bild mit zoom
 		\def\axisarguments{
 			name=A2,
 			at={(A1.south)},
 			anchor=north,
 			yshift=-1cm,
 			xmin=\ausschnittxmin,
 			xmax=\ausschnittxmax,
 			ymin=\ausschnittymin,
 			ymax=\ausschnittymax,			
 			}
 		\def\postcode{\legend{}}	
 		\input{Tikzbilder/Splitcircle-Schock/V1-plot-H.tex}		
 		% % drittes bild Q daten
 		\def\axisarguments{
 			name=A3,
 			at={(A2.south)},
 			anchor=north,
 			yshift=-1cm,	
 		}
 		\input{Tikzbilder/Splitcircle-Schock/V1-plot-Q.tex}			
 	    \end{tikzpicture}
 	}
 \caption{Schock investigation \ref{sec:capturing shock waves}: Solution of the ODE in $V_1$ over time.}
 \label{fig:V1 ODE}	    
 \end{figure}

 The evolution of the states in the ODE of vertex $V_1$ is depicted in figure \figref{fig:V1 ODE}. 
 When zooming in, 
 we can see that the $6$-th order scheme on the coarse grid of $50$ cells is closer to the reference solution computed on a grid of $200$ cells
 compared to the first order scheme. 
 Even though the initial data is non smooth,
 which reduces the order of convergence to one, the solution benefits from the high order treatment.
%  We conclude that just as for PDEs high order methods pay off even if the solution is not smooth. 
 
 In \figref{fig:E12 PDE}  the solutions along edges $E_1$ and $E_2$ are shown.
 We observe that the shock emerging from the vertex profits greatly in sharpness from a higher order scheme as well. 
 At $t=4$ the shock emerging from $V_2$ into $E_4$ after passing through $E_1$ and $E_2$ is shown in \figref{fig:E4 PDE}. 
 Again the solution of the $6$-th order scheme is much better than its first order counterpart.

 \begin{figure}
 	\centering
 	\externaltikz{Splitcircle-shock-E1-t-2}{
 	\def\ausschnittxmin{13}
 	\def\ausschnittxmax{18}
 	\def\ausschnittymin{5}
 	\def\ausschnittymax{5.7}	
 	\begin{tikzpicture}
 		\def\axisarguments{name=A1}
 		\def\precode{
 			\draw[thick,fill=gray!10,style=dotted] (axis cs: \ausschnittxmin,\ausschnittymin) rectangle (axis cs: \ausschnittxmax,\ausschnittymax);
 		}
 		\def\postcode{}	
 		\input{Tikzbilder/Splitcircle-Schock/E1-plot-h.tex}		
 		% % zweites bild mit zoom
 		\def\axisarguments{
 			name=A2,
 			at={(A1.south)},
 			anchor=north,
 			yshift=-1cm,
 			xmin=\ausschnittxmin,
 			xmax=\ausschnittxmax,
 			ymin=\ausschnittymin,
 			ymax=\ausschnittymax,			
 		}
 		\def\postcode{\legend{}}	
 		\input{Tikzbilder/Splitcircle-Schock/E1-plot-h.tex}	
 		% % drittes bild Q daten
 		\def\axisarguments{
 			name=A3,
 			at={(A2.south)},
 			anchor=north,
 			yshift=-1cm,	
 		}
 		\def\precode{}
 		\input{Tikzbilder/Splitcircle-Schock/E1-plot-q.tex}				
 	\end{tikzpicture}
 	}	
 %\caption{Solution of the PDE on $E_1$/$E_2$ at $t=2$}
  \caption{Schock investigation \ref{sec:capturing shock waves}: Solution of the PDE on   $E_1$/$E_2$ at $t=2$.}
 \label{fig:E12 PDE}
 \end{figure}

 \begin{figure}
 	\centering
 	\externaltikz{Splitcircle-shock-E4-t-4}{
 	\def\ausschnittxmin{2}
 	\def\ausschnittxmax{10}
 	\def\ausschnittymin{4.9}
 	\def\ausschnittymax{6}	
 	\begin{tikzpicture}
 	\def\pfad{Daten/SchockODE/E4/}
 	\def\axisarguments{name=A1}
 	\def\precode{
 		\draw[thick,fill=gray!10,style=dotted] (axis cs: \ausschnittxmin,\ausschnittymin) rectangle (axis cs: \ausschnittxmax,\ausschnittymax);
 	}
 	\def\postcode{}	
 	\input{Tikzbilder/Splitcircle-Schock/E4-plot-h.tex}			
 	% % zweites bild mit zoom
 	\def\axisarguments{
 		name=A2,
 		at={(A1.south)},
 		anchor=north,
 		yshift=-1cm,
 		xmin=\ausschnittxmin,
 		xmax=\ausschnittxmax,
 		ymin=\ausschnittymin,
 		ymax=\ausschnittymax,			
 	}
 	\def\postcode{\legend{}}	
 	\input{Tikzbilder/Splitcircle-Schock/E4-plot-h.tex}	
 	% % drittes bild Q daten
 	\def\axisarguments{
 		name=A3,
 		at={(A2.south)},
 		anchor=north,
 		yshift=-1cm,	
 	}
 	\def\precode{}
 	\input{Tikzbilder/Splitcircle-Schock/E4-plot-q.tex}			
 	\end{tikzpicture}
 	}	
% 	\caption{Solution of the PDE on $E_4$ at $t=4$}
	\caption{Schock investigation \ref{sec:capturing shock waves}: Solution of the PDE on $E_4$ at $t=4$.}
 	\label{fig:E4 PDE}
 \end{figure}
 \MyFloatBarrier
 
 \subsection{Shocks and lumped parameter models}
 \label{sec:Shocks LPM}
 In this test we consider a larger network of $32$ edges and $24$ nodes.
 \begin{figure}
 	\centering
 	\externaltikz{Split-and-join-LPM-network-large}{
 		\def\innereKante{blue}	
 		\begin{tikzpicture}[xscale=0.7]
 		\tikzstyle{knoten}=[circle, fill=gray!30, scale=0.6];
 		\tikzstyle{kante}=[-,line width=1];
 		\node[knoten] (1) at (-2,0) {$V_{1}$};
 		\node[knoten] (2) at (0,0) {$V_{2}$};
 		\node[knoten] (3) at (2,-2.2222) {$V_{3}$};
 		\node[knoten] (4) at (2,2.2222) {$V_{4}$};
 		\node[knoten] (5) at (4,-3.3333) {$V_{5}$};
 		\node[knoten] (6) at (4,-1.1111) {$V_{6}$};
 		\node[knoten] (7) at (4,1.1111) {$V_{7}$};
 		\node[knoten] (8) at (4,3.3333) {$V_{8}$};
 		\node[knoten] (9) at (6,-3.8889) {$V_{9}$};
 		\node[knoten] (10) at (6,-2.7778) {$V_{10}$};
 		\node[knoten] (11) at (6,-1.6667) {$V_{11}$};
 		\node[knoten] (12) at (6,-0.55556) {$V_{12}$};
 		\node[knoten] (13) at (6,0.55556) {$V_{13}$};
 		\node[knoten] (14) at (6,1.6667) {$V_{14}$};
 		\node[knoten] (15) at (6,2.7778) {$V_{15}$};
 		\node[knoten] (16) at (6,3.8889) {$V_{16}$};
 		\node[knoten] (17) at (8,-3.3333) {$V_{17}$};
 		\node[knoten] (18) at (8,-1.1111) {$V_{18}$};
 		\node[knoten] (19) at (8,1.1111) {$V_{19}$};
 		\node[knoten] (20) at (8,3.3333) {$V_{20}$};
 		\node[knoten] (21) at (10,-2.2222) {$V_{21}$};
 		\node[knoten] (22) at (10,2.2222) {$V_{22}$};
 		\node[knoten] (23) at (12,0) {$V_{23}$};
 		\node[knoten] (24) at (14,0) {$V_{24}$};
 		\draw[kante] (1) .. controls(-1.3333,0)and(-0.66667,0).. (2) node[pos=0.5,above]{$E_{1}$};
 		\draw[kante] (1) .. controls(-4.0045,-1.6032)and(-3.6622,1.5818).. (1) node[pos=0.5,above]{$E_{2}$};
 		\draw[kante] (2) .. controls(0.66667,-0.74074)and(1.3333,-1.4815).. (3) node[pos=0.5,above]{$E_{3}$};
 		\draw[kante] (2) .. controls(0.66667,0.74074)and(1.3333,1.4815).. (4) node[pos=0.5,above]{$E_{4}$};
 		\draw[kante,\innereKante] (3) .. controls(2.6667,-2.5926)and(3.3333,-2.963).. (5) node[pos=0.5,above]{$E_{5}$};
 		\draw[kante,\innereKante] (3) .. controls(2.6667,-1.8519)and(3.3333,-1.4815).. (6) node[pos=0.5,above]{$E_{6}$};
 		\draw[kante,\innereKante] (4) .. controls(2.6667,1.8519)and(3.3333,1.4815).. (7) node[pos=0.5,above]{$E_{7}$};
 		\draw[kante,\innereKante] (4) .. controls(2.6667,2.5926)and(3.3333,2.963).. (8) node[pos=0.5,above]{$E_{8}$};
 		\draw[kante,\innereKante] (5) .. controls(4.6667,-3.5185)and(5.3333,-3.7037).. (9) node[pos=0.5,above]{$E_{9}$};
 		\draw[kante,\innereKante] (5) .. controls(4.6667,-3.1481)and(5.3333,-2.963).. (10) node[pos=0.5,above]{$E_{10}$};
 		\draw[kante,\innereKante] (6) .. controls(4.6667,-1.2963)and(5.3333,-1.4815).. (11) node[pos=0.5,above]{$E_{11}$};
 		\draw[kante,\innereKante] (6) .. controls(4.6667,-0.92593)and(5.3333,-0.74074).. (12) node[pos=0.5,above]{$E_{12}$};
 		\draw[kante,\innereKante] (7) .. controls(4.6667,0.92593)and(5.3333,0.74074).. (13) node[pos=0.5,above]{$E_{13}$};
 		\draw[kante,\innereKante] (7) .. controls(4.6667,1.2963)and(5.3333,1.4815).. (14) node[pos=0.5,above]{$E_{14}$};
 		\draw[kante,\innereKante] (8) .. controls(4.6667,3.1481)and(5.3333,2.963).. (15) node[pos=0.5,above]{$E_{15}$};
 		\draw[kante,\innereKante] (8) .. controls(4.6667,3.5185)and(5.3333,3.7037).. (16) node[pos=0.5,above]{$E_{16}$};
 		\draw[kante,\innereKante] (9) .. controls(6.6667,-3.7037)and(7.3333,-3.5185).. (17) node[pos=0.5,above]{$E_{17}$};
 		\draw[kante,\innereKante] (10) .. controls(6.6667,-2.963)and(7.3333,-3.1481).. (17) node[pos=0.5,above]{$E_{18}$};
 		\draw[kante,\innereKante] (11) .. controls(6.6667,-1.4815)and(7.3333,-1.2963).. (18) node[pos=0.5,above]{$E_{19}$};
 		\draw[kante,\innereKante] (12) .. controls(6.6667,-0.74074)and(7.3333,-0.92593).. (18) node[pos=0.5,above]{$E_{20}$};
 		\draw[kante,\innereKante] (13) .. controls(6.6667,0.74074)and(7.3333,0.92593).. (19) node[pos=0.5,above]{$E_{21}$};
 		\draw[kante,\innereKante] (14) .. controls(6.6667,1.4815)and(7.3333,1.2963).. (19) node[pos=0.5,above]{$E_{22}$};
 		\draw[kante,\innereKante] (15) .. controls(6.6667,2.963)and(7.3333,3.1481).. (20) node[pos=0.5,above]{$E_{23}$};
 		\draw[kante,\innereKante] (16) .. controls(6.6667,3.7037)and(7.3333,3.5185).. (20) node[pos=0.5,above]{$E_{24}$};
 		\draw[kante,\innereKante] (17) .. controls(8.6667,-2.963)and(9.3333,-2.5926).. (21) node[pos=0.5,above]{$E_{25}$};
 		\draw[kante,\innereKante] (18) .. controls(8.6667,-1.4815)and(9.3333,-1.8519).. (21) node[pos=0.5,above]{$E_{26}$};
 		\draw[kante,\innereKante] (19) .. controls(8.6667,1.4815)and(9.3333,1.8519).. (22) node[pos=0.5,above]{$E_{27}$};
 		\draw[kante,\innereKante] (20) .. controls(8.6667,2.963)and(9.3333,2.5926).. (22) node[pos=0.5,above]{$E_{28}$};
 		\draw[kante] (21) .. controls(10.6667,-1.4815)and(11.3333,-0.74074).. (23) node[pos=0.5,above]{$E_{29}$};
 		\draw[kante] (22) .. controls(10.6667,1.4815)and(11.3333,0.74074).. (23) node[pos=0.5,above]{$E_{30}$};
 		\draw[kante] (23) .. controls(12.6667,0)and(13.3333,0).. (24) node[pos=0.5,above]{$E_{31}$};
 		\draw[kante] (24) .. controls(15.49,-1.8699)and(15.387,1.4273).. (24) node[pos=0.5,above]{$E_{32}$};
 		\draw[fill=green,opacity=0.2] (6,-2.2) ellipse (5 and 2.2);
 		\end{tikzpicture}
 	}
 	\caption{Split and join network with lumping of the green area.}
 	\label{fig:Split and join network large with lumping}
 \end{figure}
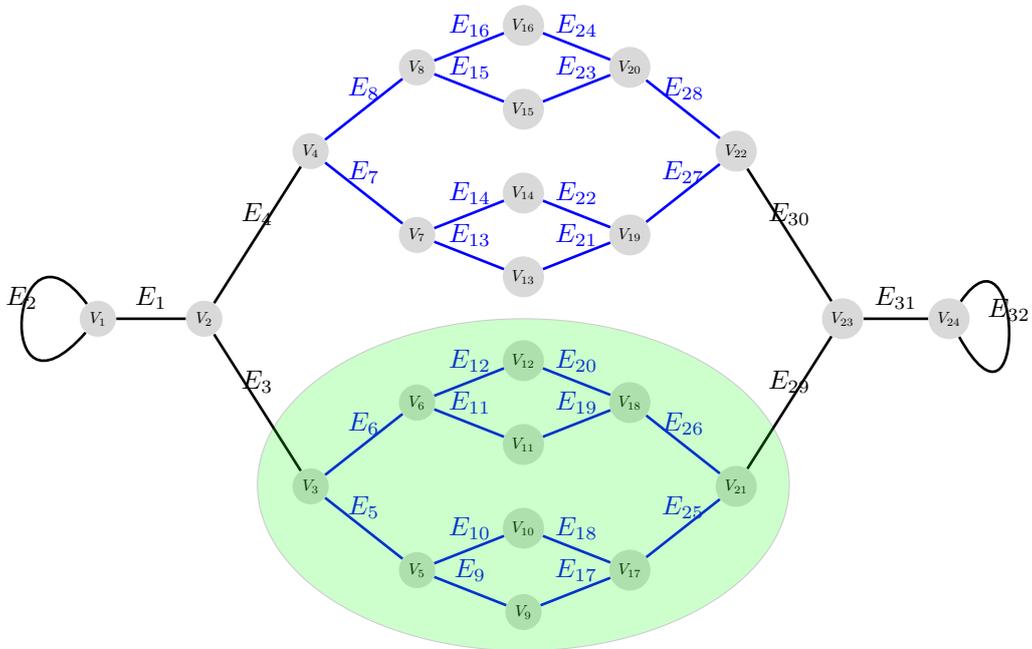
 As shown in \figref{fig:Split and join network large with lumping}, the network is of a tree like structure,
 inspired by the human circulatory system e.g. \cite{muller2014global}.
 The first four edges $E_1,E_2,E_3,E_4$ and the last four $E_{29},E_{30},E_{31},E_{32}$ have a length of $L=25$,
 for all remaining edges we choose $L=2.5$.
 In the nodes equal height coupling conditions \eqref{eq:equalheights} are used.
 In order to investigate the influence of parameter lumping on the solution 
 we model the lower half of the network by an ODE as described in Section \ref{S:Lumped Parameter Models}.

 As initial conditions we choose water at rest on all edges.
 On $E_1$ we impose the following Riemann Problem 
 \begin{align*}
 h^1(0,x)&=
 \begin{cases}
 3 & x<18.5 \\
 2 & \text{else}
 \end{cases}
 &q^1(0,x)&=0\ ,
 \end{align*}
whereas the remaining edges continue these constant states
 \begin{align*}
 u^2(0,x)&
 \equiv (3,0)^T\ ,
%  \begin{pmatrix}
%  3 \\ 0
%  \end{pmatrix}
 & u^i(0,x)&\equiv
  (2,0)^T
%   \begin{pmatrix}
%  2 \\ 0
%  \end{pmatrix}
\quad \forall i\geq 3\ .
 \end{align*}
  Due to the symmetric nature of this setting and without applying the lumping to the lower part of the network,
  the solution in some branches of the network coincide .
  In this case we have the following identities
 \begin{align*}
 u^3&=u^4 \ , \\
 u^5=u^6&=u^7=u^8 \ , \\
 u^9=u^{10}=u^{11}=u^{12}&= u^{13}=u^{14}=u^{15}=u^{16}\ , \\
 u^{17}=u^{18}=u^{19}= u^{20}&=u^{21}=u^{22}=u^{23}=u^{24}\ , \\
 u^{25}=u^{26}&=u^{27}=u^{28}\ , \\
 u^{29}&=u^{30}\ .
 \end{align*}
 Therefore it suffices to look at one edge of each group
 and we can directly compare the solutions of the lumped part of the network to its PDE counterparts.
 \begin{figure}
 	\centering
 	\externaltikz{LPM-height-vs-full-PDE}{
 		\begin{tikzpicture}
 		\def\seta{dashed}
 		\def\setb{smooth}
 		\def\axisarguments{height=4cm,yshift=-0.3cm,xmax=15,ymax=2.4}
 		\def\pfad{Daten/LPMHalskette/LPM-vs-PDE/}
 		\def\nameanfang{}
 		\def\nameende{.txt}
 		\begin{axis}[
 		%	ylabel=y axis label,
 		width=\textwidth,
 		xmin=5,
 		xmax=10,
 		name=A1,
 		xticklabels={,,},
 		legend style={at={(0.5,1)},xshift=-0cm,yshift=0.2cm,anchor=center,nodes=right},
 		legend columns=2,
 		\axisarguments
 		]% file
 		\addplot[\seta,no marks] plot file {\pfad Dichte_mittel_E8\nameende};
 		\addlegendentry{Average Height PDE}
 		\addplot[\setb,no marks] plot file {\pfad LPM-Data-cells-100-order-6-component1\nameende};
 		\addlegendentry{Height LPM}	
 		\end{axis}
 		\begin{axis}[
 		width=\textwidth,
 		xmin=5,
 		xmax=10,
 		name=A2,
 		at={(A1.south)},
 		anchor=north,
 		xticklabels={,,},
 		\axisarguments
 		]
 		\addplot[\seta,no marks] plot file {\pfad Dichte_mittel_E16\nameende};
 		\addplot[\setb,no marks] plot file {\pfad LPM-Data-cells-100-order-6-component5\nameende};
 		\end{axis}
 		\begin{axis}[
 		width=\textwidth,
 		xmin=5,
 		xmax=10,
 		name=A3,
 		at={(A2.south)},
 		anchor=north,
 		xticklabels={,,},
 		\axisarguments
 		]
 		\addplot[\seta,no marks] plot file {\pfad Dichte_mittel_E24\nameende};
 		\addplot[\setb,no marks] plot file {\pfad LPM-Data-cells-100-order-6-component13\nameende};
 		\end{axis}
 		\begin{axis}[
 		xlabel=$t$,
 		width=\textwidth,
 		xmin=5,
 		xmax=10,
 		name=A4,
 		at={(A3.south)},
 		anchor=north,
 		\axisarguments
 		]
 		\addplot[\seta,no marks] plot file {\pfad Dichte_mittel_E28\nameende};	
 		\addplot[\setb,no marks] plot file {\pfad LPM-Data-cells-100-order-6-component21\nameende};
 		\end{axis}
 		\end{tikzpicture}
 	}
 	\caption{Shocks and LPM \ref{sec:Shocks LPM}: Height components of the LPM model over time compared to the averaged solution of the full PDE simulation. From top to bottom: $E_{8}$, $E_{16}$, $E_{24}$ and $E_{28}$.}
 %	\caption{Height components of the LPM model over time compared to the averaged solution of the full PDE simulation. From top to bottom: $E_{8}$, $E_{16}$, $E_{24}$ and $E_{28}$}
 	\label{fig:LPM height vs full PDE simulation}
 \end{figure}
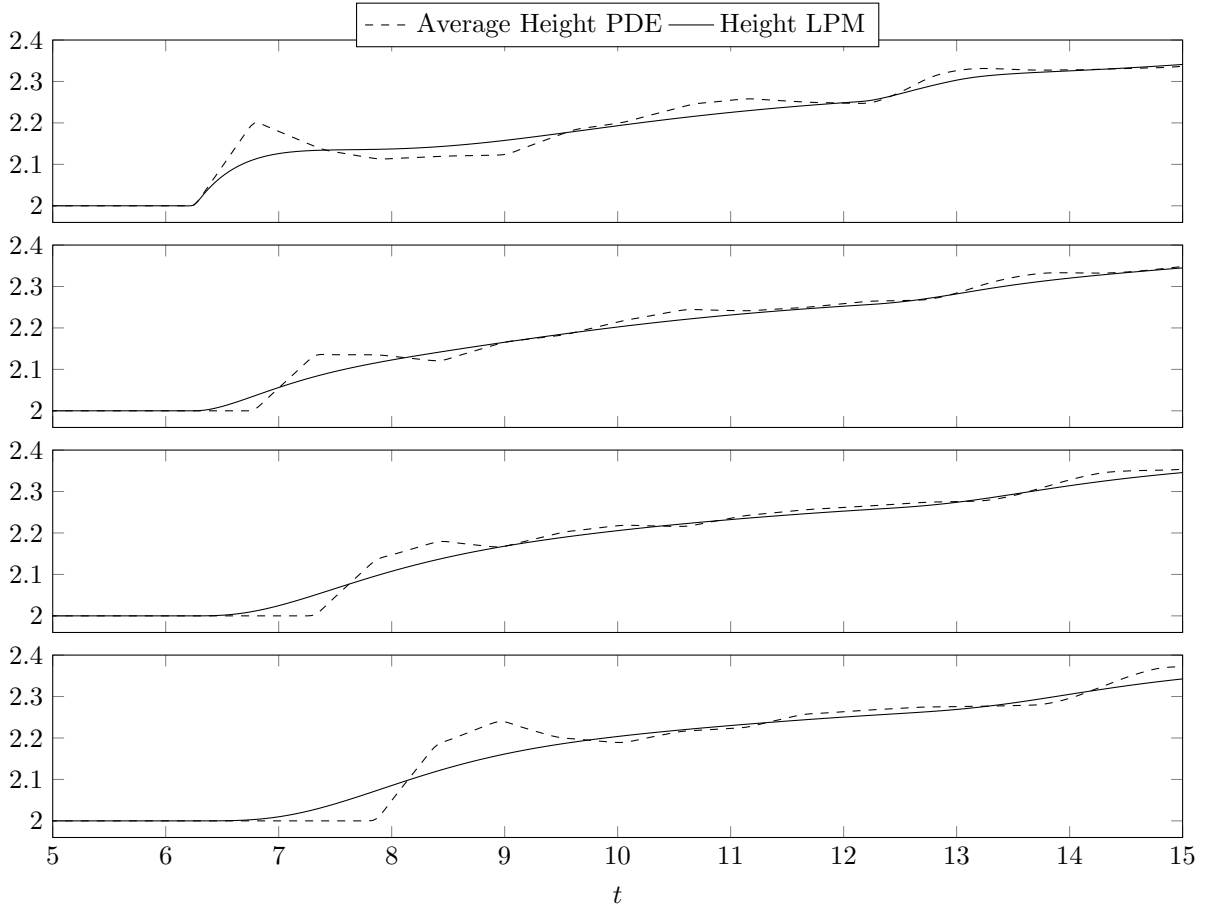
 
 \begin{figure}
 	\centering
 	\externaltikz{LPM-impulse-vs-full-PDE}{
 		\begin{tikzpicture}
 		\def\seta{dashed}
 		\def\setb{smooth}
 		\def\axisarguments{height=4cm,yshift=-0.3cm,xmax=15}
 		\def\pfad{Daten/LPMHalskette/LPM-vs-PDE/}
 		\def\nameanfang{}
 		\def\nameende{.txt}
 		\begin{axis}[
 		%	ylabel=y axis label,
 		width=\textwidth,
 		xmin=5,
 		xmax=10,
 		name=A1,
 		xticklabels={,,},
 		legend style={at={(0.5,1)},xshift=-0cm,yshift=0.2cm,anchor=center,nodes=right},
 		legend columns=2,
 		\axisarguments
 		]% file
 		\addplot[\seta,no marks] plot file {\pfad Impulse_mittel_E8\nameende};
 		\addlegendentry{Average Impulse PDE}
 		\addplot[\setb,no marks] plot file {\pfad LPM-Data-cells-100-order-6-component2\nameende};
 		\addlegendentry{Impulse LPM}	
 		\end{axis}
 		\begin{axis}[
 		width=\textwidth,
 		xmin=5,
 		xmax=10,
 		name=A2,
 		at={(A1.south)},
 		anchor=north,
 		xticklabels={,,},
 		\axisarguments
 		]
 		\addplot[\seta,no marks] plot file {\pfad Impulse_mittel_E16\nameende};
 		\addplot[\setb,no marks] plot file {\pfad LPM-Data-cells-100-order-6-component6\nameende};
 		\end{axis}
 		\begin{axis}[
 		width=\textwidth,
 		xmin=5,
 		xmax=10,
 		name=A3,
 		at={(A2.south)},
 		anchor=north,
 		xticklabels={,,},
 		\axisarguments
 		]
 		\addplot[\seta,no marks] plot file {\pfad Impulse_mittel_E24\nameende};
 		\addplot[\setb,no marks] plot file {\pfad LPM-Data-cells-100-order-6-component14\nameende};
 		\end{axis}
 		\begin{axis}[
 		xlabel=$t$,
 		width=\textwidth,
 		xmin=5,
 		xmax=10,
 		name=A4,
 		at={(A3.south)},
 		anchor=north,
 		\axisarguments
 		]
 		\addplot[\seta,no marks] plot file {\pfad Impulse_mittel_E28\nameende};	
 		\addplot[\setb,no marks] plot file {\pfad LPM-Data-cells-100-order-6-component22\nameende};
 		\end{axis}
 		\end{tikzpicture}
 	}
 	\caption{Shocks and LPM \ref{sec:Shocks LPM}: Impulse components of the LPM model over time compared to the averaged solution of the full PDE simulation. From top to bottom: $E_{8}$, $E_{16}$, $E_{24}$ and $E_{28}$.} 	
% 	\caption{Impulse components of the LPM model over time compared to the averaged solution of the full PDE simulation. From top to bottom: $E_{8}$, $E_{16}$, $E_{24}$ and $E_{28}$}
 	\label{fig:LPM impulse vs full PDE simulation}
 \end{figure}
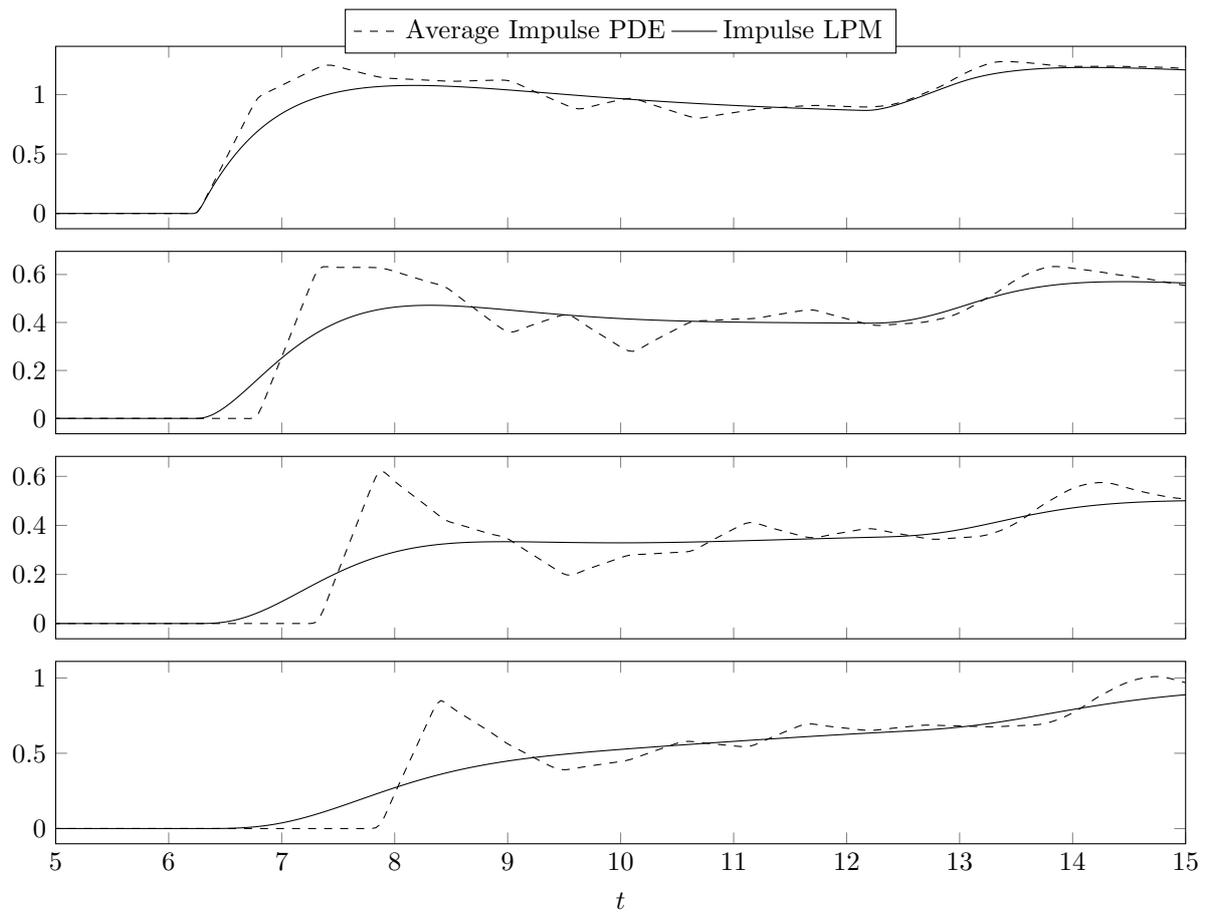
 In \figref{fig:LPM height vs full PDE simulation} we show the heights $H$ in the LPM vertex and 
 the averaged heights on a corresponding edge located on the upper half of the network.
 Analogously the momentum components can be seen in \figref{fig:LPM impulse vs full PDE simulation}. 
 
 Due to the very coarse spatial resolution the states of the LPM model can not resolve the incoming shock wave accurately. 
 Despite this strong diffusion caused by the model, the ODE captures the general behavior of the flow.
 \begin{figure}
 	\centering
 	\externaltikz{LPM-Halskette-solution-after-LPM}{
 	\begin{tikzpicture}% file
 	\def\pfad{Daten/LPMHalskette/Edge-after-LPM/}
 	\begin{axis}[
 	xlabel=$x$,
 	ylabel=$h$,
 	xmin=0,
 	xmax=25,
 	width=0.45\textwidth,
 	name=A1,
 	]
 	\addplot[black] file {\pfad Edge-After-LPM-cells-50-order-1-h.txt};
 	\addlegendentry{order 1}
 	\addplot[blue] file {\pfad Edge-After-LPM-cells-50-order-6-h.txt};
 	\addlegendentry{order 6}
 	\addplot[red] file {\pfad Edge-After-LPM-reference-cells-100-order-6-h.txt};
 	\addlegendentry{Reference}
 	\end{axis}
 	\begin{axis}[
 	xlabel=$x$,
 	ylabel=$q$,
 	xmin=0,
 	xmax=25,
 	width=0.45\textwidth,
 	name=A1,
 	at={(A1.east)},
 	anchor=west,
 	xshift=2cm,
 	]
 	\addplot[black] file {\pfad Edge-After-LPM-cells-50-order-1-q.txt};
 	\addlegendentry{order 1}
 	\addplot[blue] file {\pfad Edge-After-LPM-cells-50-order-6-q.txt};
 	\addlegendentry{order 6}
 	\addplot[red] file {\pfad Edge-After-LPM-reference-cells-100-order-6-q.txt};
 	\addlegendentry{Reference}
 	\end{axis}
 	\end{tikzpicture}
 	}
% 	\caption{Solution on $E_{29}$ after the LPM vertex compared to the reference solution on $E_{30}$}
 	\caption{Shocks and LPM \ref{sec:Shocks LPM}: Solution on $E_{29}$ after the LPM vertex compared to the reference solution on $E_{30}$.}
 	
 	\label{fig:LPM Halskette solution after LPM vertex}
 \end{figure}
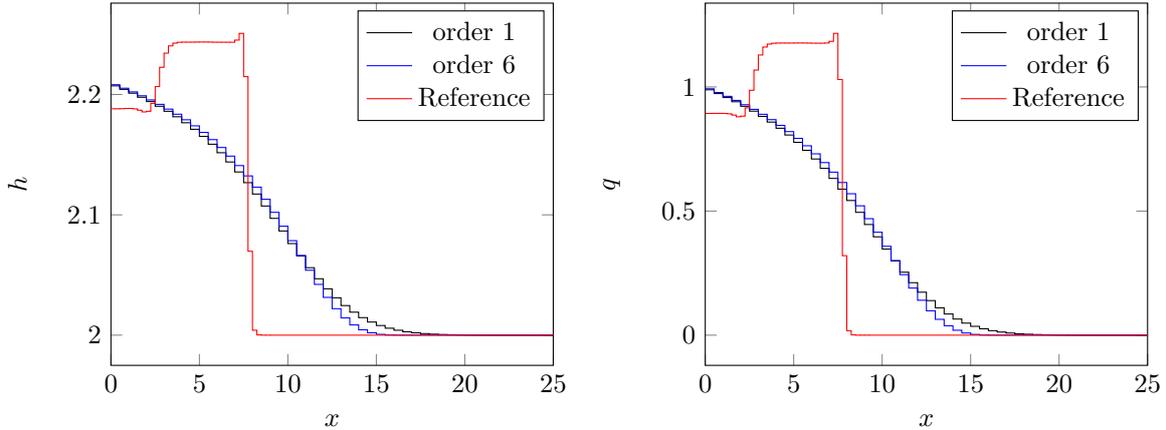
 
 The correct capturing of the wave speed can be observed in \figref{fig:LPM Halskette solution after LPM vertex}.
 Here the solution on the edges $E_{29}$ and $E_{30}$ are shown. 
 This provides a comparison between the shock wave, that has passed the lumped branch of the network
 and wave transported by the PDE model.
 Clearly the solution of the PDE model can resolve much finer structures.
 
 \MyFloatBarrier

\subsubsection{Sources in Lumped Parameter Models}
\label{sec:Soures and well balancing}
\begin{figure}
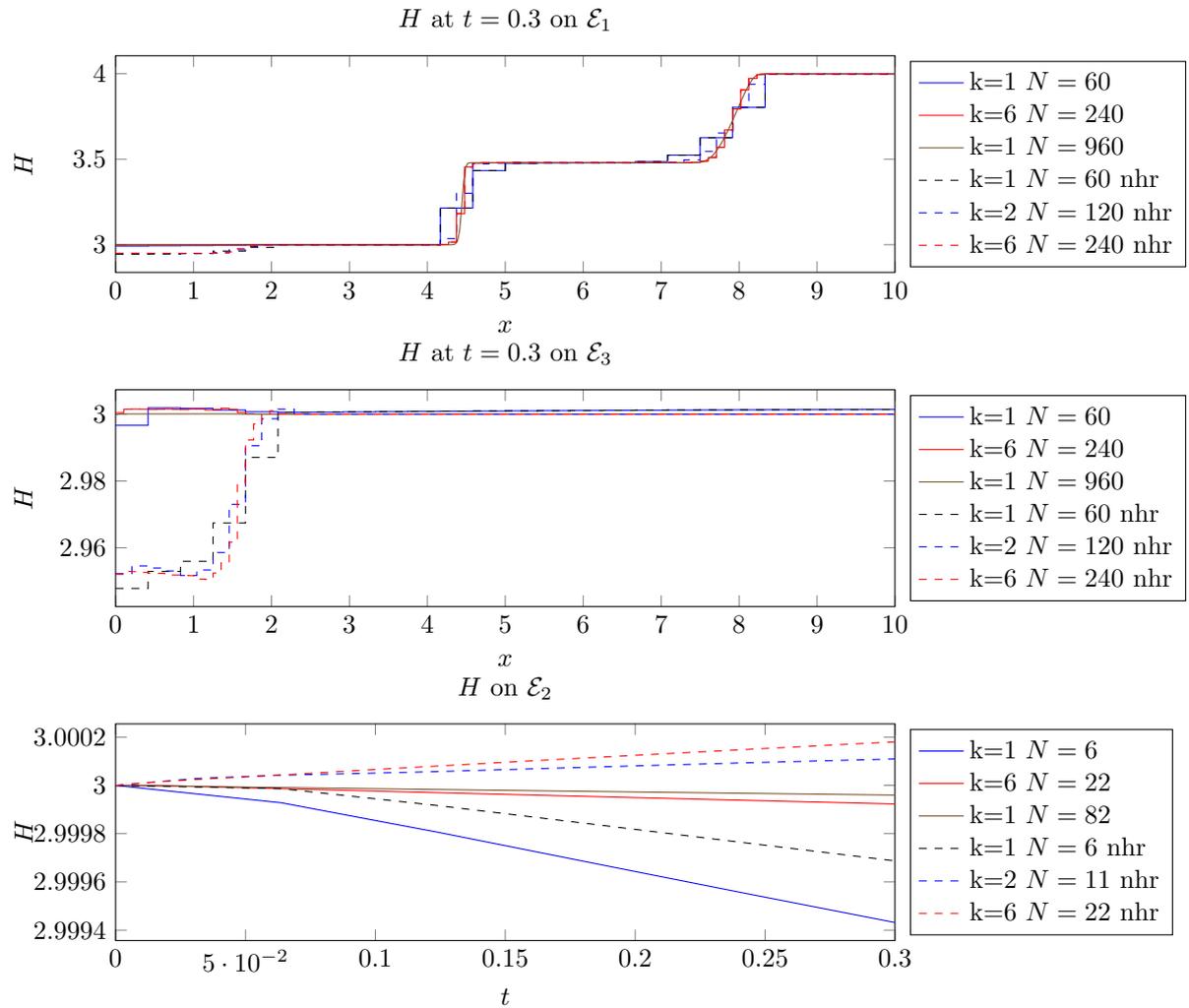

 	\centering
 	\externaltikz{LPM-well-balanced-example-H}{
 		\relinput{Tikzbilder/NUM_splitcircle_bottom_LPM/H-bild-diss/H-bild-diss.tex}		
 	}
 	%\caption{Kante 1}
 	\caption{Well balanced LPM \ref{sec:Soures and well balancing}: $H$ on $E_1$.}
 	\label{fig:LPM-wb-H}
 \end{figure}

To demonstrate the necessity of hydrostatic reconstruction in LPM models,
 we simulate the following initial value problem consisting of shallow water equations on a split circle network with bottom elevation. $E_2$ and the two vertices are turned into a LPM. We use bottom profiles of linear, polynomial and trigonometrical type:
\begin{align*}
b_1(x)&=0.3\frac{x}{25} 
& 
b_2(x)&=0.3\left(
\frac{x}{25}+
\left(
\frac{x}{25}-\frac{1}{2}
\right)^2-\frac{1}{4}
\right)
&
b_3(x)&=
0.3\left(
\frac{x}{25}+\sin(\pi \frac{x}{25})
\right)
\\
H_1(x)&=3+\mathbf{1}_{[\frac{1}{4}25 , \frac{3}{4}25]}(x)
&H_2(x)&\equiv 3
&H_3(x)&\equiv 3
\\
\rho_1(x)&=H_1(x)-b_1(x) 
&\rho_2(x)&=H_2(x)-b_2(x) 
&\rho_3(x)&=H_3(x)-b_3(x) 
\end{align*}
The results of the simulations are shown in figure \ref{fig:LPM-wb-H} and \ref{fig:LPM-wb-Q}. The solutions on the PDE edges are shown at $t=0.3$, the state of the ODE over the complete time interval. In the interior of the domain we can observe the asymptotically well balancedness of ADER schemes as reported in \cite{castro2008solvers}, the slightly stronger oscillations at the vertices are caused by the higher sensitivity of the coupling conditions amplifying the unavoidable errors stemming from the one sided polynomial reconstruction. The big oscillations however are avoided as expected.   
%Note the oscillations at the boundaries, cause: errors in reconstruction in concert with the higher sensitivity of the coupling

 %ggf rausnehmen
 \begin{figure}
 	\centering
 	\externaltikz{LPM-well-balanced-example-Q}{
 		\relinput{Tikzbilder/NUM_splitcircle_bottom_LPM/Q-bild-diss/Q-bild-diss.tex}		
 	}
% 	\caption{Kante 1}
 	\caption{Well balanced LPM \ref{sec:Soures and well balancing}: $q$ on $E_1$.} 	
 	\label{fig:LPM-wb-Q}
 \end{figure}

% \begin{figure}
% 	\centering
% 	\externaltikz{LPM-well-balanced-example-E1}{
% 		\relinput{Tikzbilder/NUM_splitcircle_bottom_LPM/Kante1_PDE/Kante1_PDE.tex}		
% 	}
% 	\caption{Kante 1}
% 	\label{fig:LPM-wb-E1}
% \end{figure}
%
% \begin{figure}
% 	\centering
% 	\externaltikz{LPM-well-balanced-example-E3}{
% 		\relinput{Tikzbilder/NUM_splitcircle_bottom_LPM/Kante3_PDE/Kante3_PDE.tex}		
% 	}
% 	\caption{Kante 3}
% 	\label{fig:LPM-wb-E3}
% \end{figure}
%
%  \begin{figure}
%  	\centering
%  	\externaltikz{LPM-well-balanced-example-E2}{
%  		\relinput{Tikzbilder/NUM_splitcircle_bottom_LPM/Kante2_LPM/Kante2_LPM.tex}		
%  	}
%  	\caption{Kante 2 LPM}
%  	\label{fig:LPM-wb-E2}
%  \end{figure}

%%%%%%%%%%%%%%%%%%%

\section{Conclusion}
\label{S:Conclusion}
High order GRP solvers for ODE vertices and LPM models of \TT and \HEOC type are introduced in this work. Extensive tests showed that they indeed exhibit the high order of convergence they were designed for. Numerical examples showed that these technique can indeed be used to build very accurate and stable numerical methods for networks of conservation laws including vertices with ODEs and lumped parameter models. Especially the solver of \HEOC type is a vast improvement in terms of applicability over the \TT approach introduced in \cite{borsche2014ader}.  
	

\bibliographystyle{elsarticle-num}
\bibliography{Qiqqa2BibTexExport,paper_ODE_coupling}


\end{document}